\documentclass[a4paper,10pt]{article}
\usepackage{amssymb,graphicx,graphics,amsmath}
\usepackage{makeidx}
\makeindex
\newtheorem{next}{}[section]
\newcommand{\nx}[1]{{\bf #1.}\ }
\newcommand{\proof}{\indent{\it Proof.}\quad}
\newcommand{\qed}{\hfill\rule{1.6mm}{1.6mm}\,}%{\hfill$\Box$}
\newcommand{\qeds}{\begin{flushright}
$\Box$ 
\end{flushright}}%{\hfill$\Box$}

\newcommand{\B}{{\cal B}}
\newcommand {\D} {\Delta}
\renewcommand {\k} {\Bbbk}

\newcommand {\Z} {{\mathbb Z}}
\newcommand {\N} {{\mathbb N}}

\newcommand {\FF}{{\cal F}}

\newcommand{\si}{\k[{\bf x}]}
\newcommand{\sis}{\k[{\bf z}]}

\newcommand{\M}{{\cal M}}
\newcommand{\Dd}{{\cal D}}
\newcommand{\I}{{\cal I}}
\newcommand{\J}{{\cal J}}

\newcommand{\Cc}{{\mathfrak C}}
\newcommand{\Qq}{{\cal Q}}
\newcommand{\Pp}{{\cal P}}

\newcommand{\dtb}{$\stackrel{\mbox{\circle*{4.0}}}{}$}

\newcommand{\igno} [1] {}

\newcommand{\atrib}[1]{$(${\rm #1}$)$}

\title{The ${\rm N}_{2,p}$-property of binomial extensions of simplicial complexes}
\author{Hernan de Alba Casillas, Marcel Morales}
\date{}

\begin{document}
\maketitle

{\small\abstract{M. Morales introduced a family of binomial ideals that are binomial extensions of square free monomial ideals. Let $I\subset \si$ be a square free monomial ideal and $J\subset\sis$ a sum of scroll ideals with some extra conditions, we define the binomial extension of $I$ as $\B=I+J\subset \sis$. We set $p_2(\B)$ the minimal $i\in\N$ such that there exists $j>2$ such that $\beta_{i,i+j}(\B)\neq 0$. In the case where $J=0$, Fr\"oberg characterized combinatorally the case $p_2(I)=\infty$; later Eisenbud et al. solved the case $p_2(I)<\infty$. We obtain a similar result as Fr\"oberg for the binomial extensions and we find lower and upper bounds of $p_2(\B)$ for some families of binomial extensions in combinatorial terms as Eisenbud et al. With some additional hypothesis we can compute $p_2(\B)$.}}

%------------------------
%\section{Introduction}

%------------------------
\section{Preliminaries}
%--------------------
\subsection*{${\bf N}_{d,p}$-property}

Let $V$ be a set of variables and $\si=\k[V]$ the polynomial ring on the set of variables $V$ over a field $\k$, multigraded by a commutative cancellative semigroup $\Sigma$ with unity and without inversible elements excepts the unity. We simply say $\si$ is $\Sigma$-(positively) multigraded. Let $I\subset\si$ a $\Sigma$-multigraded ideal of $\si$. A $\Sigma$-multigraded minimal free resolution of $I$ is given by:
$$0\rightarrow \bigoplus_{\begin{array}{c}j=1\\ a_{j,\rho}\in \Sigma \end{array}}^{n_\rho} \si(-a_{j,\rho})^{\beta_{\rho,a_{j,\rho}}}\rightarrow \dots\rightarrow \bigoplus_{\begin{array}{c}j=1\\ a_{j,0}\in \Sigma \end{array}}^{n_0} \si(-a_{j,0})^{\beta_{0,a_{j,0}}}\rightarrow I\rightarrow 0,$$
where $\beta_{i,a_j}\in\N^*$ and $a_j\in \Sigma$. The ranks $\beta_{i,a_j}(I)=\beta_{i,a_j}$ are called $\Sigma$-multigraded (or simply $\Sigma$-graded) Betti numbers. In case where $\Sigma=\N$ and $\si$ have the $\N$-standard graduation the $\N$-graded Betti numbers of $I$ will just be called the Betti-numbers of $I$.

%The following theorem is a well-known result about the Betti numbers of an ideal and the Betti numbers of its initial ideal. 

\begin{next}\nx{Theorem}\atrib{Upper-semicontinuity}\label{betti-gradue-in}
 Let $\Sigma$ be a semigroup and $I\subset\si$ be a $\Sigma$ multigraded ideal, then for all monomial order $\succ$ on $\si$ the Betti numbers of I satisfy the following inequality:
\begin{center} 
	$\beta_{i,a}(I)\leq\beta_{i,a}({\rm in}_\succ(I))$, for all $i\in \N$ and $a\in \Sigma$.
\end{center}
\end{next}
A good reference for theorem \ref{betti-gradue-in} is chapter 8 of \cite{israel}.

\begin{next}\nx{Remark}\label{betti-in-i}
{\rm If there exists a monomial order $\succ$ such that $\beta_{i,a}({\rm in}_\succ(I))=0$, then $\beta_{i,a}(I)=0$.}
\end{next}

\begin{next}\nx{Definition}\label{n2P-property}
 {\rm Let $I\subset\si$ be a homogeneous ($\N$-graded standard) ideal of $\si$ and $d,p\in\N$, $d\neq 0$ and $p\neq 0$. We say that $I$ satisfies the ${\rm N}_{d,p}$ property if and only if for all $i\leq p-1$ and $j\geq 1$ $\beta_{i,i+d+j}(I)=0$. Now, if $d\in\N$, $d\neq 0$ we set $p_d(I)\in\N$ as the maximal integer $p$ such that $I$ satisfies the ${\bf N}_{d,p}$-property, if this maximum exists, otherwise we set $p_d(I)=\infty$ and in this case we say that $I$ is $d$-regular. In case where $I$ is generated by homogeneous elements of degree $d$ and is $d$-regular, we say that $I$ is $d$-linear.}
\end{next}

\begin{next}\nx{Proposition}\label{igualdad-ini}
Let $I\subset\si$ be an ideal generated by homogeneous elements of degree $d$, $\succ$ a monomial order over $\si$, then
\begin{center}
	$\beta_{i,j}(I)=\beta_{i,j}({\rm in}(I))$ for all $i\leq p_d({\rm in}(I))-1$ and $j\in\N$; and\\
	$\beta_{p_d({\rm in}(I)),p_d({\rm in}(I))+d}(I)=\beta_{p_d({\rm in}(I)),p_d({\rm in}(I))+d}({\rm in}(I))$,
\end{center}
thus $p_d(I)\geq p_d({\rm in}(I))$. In particular, if ${\rm in}(I)$ is $d$-linear, then $I$ is also $d$-linear and\break$\beta_{i,j}(I)=\beta_{i,j}({\rm in}(I))$ for all $i,j\in\N$.
\end{next} 
\proof
We set $p=p_d({\rm in}(I))$. So the Betti-table of ${\rm in}(I)$ has the shape:
\begin{table}[h]
							\centering
								\begin{tabular}{|c|ccccccc|}\hline
										&$0$			&$1$		&$\dots$	&$p-1$&$p$&$\dots$& $\rho$\\ \hline
									$d$	&$\beta_{0,d}({\rm in}(I))$ &$\beta_{1,d+1}({\rm in}(I))$&$\dots$&$\beta_{p-1,p+d-1}({\rm in}(I))$&$\beta_{p,p+d}({\rm in}(I))$&$\dots$&$\beta_{\rho,\rho+d}({\rm in}(I))$\\ 
									$d+1$	&$0$ &$0$&$\dots$&$0$&$\beta_{p,p+d+1}({\rm in}(I))$&$\dots$&$\beta_{\rho,\rho+d+1}({\rm in}(I))$\\ 
									$\vdots$& & & & & & & 	\\
									$m$	&$0$ &$0$&$\dots$&$0$&$\beta_{p,p+m}({\rm in}(I))$&$\dots$&$\beta_{\rho,\rho+m}({\rm in}(I))$\\ \hline
								\end{tabular}\\
								\label{numerillus5}
\end{table}

By the remark \ref{betti-in-i} for all $i<p-1$ and $j\in\N$, $\beta_{i,d+i+j}(I)=0$. Let $H(\si/I;t)$ be the Hilbert series of $\si/I$, we know that  $$\frac{{\cal K}(I;t)}{(1-t)^n}=H(\si/I;t)=H(\si/{\rm in}(I);t)=\frac{{\cal K}({\rm in}(I);t)}{(1-t)^n}$$
where ${\cal K}(M;t)=\sum_{i\geq 1, \mbox{ }j\geq 0}(-1)^{i}\beta_{i-1,j}(M)t^j$ with $M=I$, or $M={\rm in}(I)$.
So ${\cal K}(I;t)={\cal K}({\rm in}(I);t)$ and consequently $[{\cal K}(I;t)]_{p+d}=[{\cal K}({\rm in}(I);t)]_{p+d}$. Hence
$$\sum_{i=0}^{p-1}(-1)^i\beta_{i,i+d}(I)t^{i+d}+\beta_{p,d+p}(I)t^{d+p}=\sum_{i=0}^{p-1}(-1)^i\beta_{i,i+d}({\rm in}(I))t^{i+d}+\beta_{p,d+p}({\rm in}(I))t^{d+p}.$$
Then, for $0\leq i\leq p$, $\beta_{i,i+d}(I)=\beta_{i,i+d}({\rm in}(I))$, and $p_d(I)\geq p_d({\rm in}(I))$.
\qed\\

\begin{next}\nx{Remark}
Let $I\subset \si$ be an ideal generated by homogeneous element of degree $\geq d$, so\break $\beta_{0,d}=\beta_{0,d}({\rm in}(I))$. Moreover, if ${\rm in}(I)$ is generated by monomials of degree $d$ then $I$ is generated by elements of degree $d$.
\end{next}
%------------------------
\subsection*{Clique complexes}

A simplicial complex $\Gamma$ over a vertex set $V(\Gamma)=\{x_1,\dots,x_n\}$ is a collection of subsets of $V(\Gamma)$ such that: for all i, the set $\{x_i\}$ is in $\Gamma$; and if $F\in \Gamma$ and $G\subset F$ then $G\in\Gamma$. An element of a simplicial complex $\Gamma$ is called a {\bf face} of $\Gamma$. \igno{The dimension of a face $F$ of $\Gamma$, denoted by ${\rm dim}\mbox{ } F$, is defined to be $|F|-1$, where $|F|$ denotes the number of vertices in $F$. The dimension of $\Gamma$, denoted by ${\rm dim}\mbox{ } \Gamma$, is defined to be the maximal dimension of a face in $\Gamma$.} The maximal faces of $\Gamma$ under inclusion are called facets of $\Gamma$. The 1-skeleton $\Gamma_1$ of $\Gamma$ is the subcomplex of $\Gamma$ consisting of all the faces of $\Gamma$ whose cardinality is $\leq 2$. Let $W\subset V(\Gamma)$, the restriction of $\Gamma$ over $W$ denoted by $\Gamma_W$ is the simplicial complex composed by all the faces $F$ of $\Gamma$ such that $F\subset W$.\igno{The support of a face $F\in\Gamma$ is the element $F=(a_1,dots,a_n)\in\N^n$ such that $a_i=0$ if $x_i\notin F$ and $a_i=1$ if $x_i\in F$.} We set $x^{F}=\prod_{x_i\in F}x_i$. The Stanley-Reisner ideal of $\Gamma$, denoted by $I_\Gamma$, is the ideal of $\k[V(\Gamma)]$ defined by $I_\Gamma=({\bf x}^{F}:F\notin \Gamma)$. 

\begin{next}\nx{Theorem}\label{Hochster's Formula}\atrib{{\bf Hochster's Formula}}
 All nonzero Betti numbers $\N^n$-graded of $I_{\Gamma}$ and $S/I_{\Gamma}$ lie in squarefree degrees $\sigma=(\sigma_1,\dots,\sigma_n)\in\N^n$, and
$$\beta_{i,\sigma}(I_{\Gamma})=\beta_{i+1,\sigma}(S/I_{\Gamma})={\rm dim}_{\k}\widetilde{H}_{|\sigma|-i-2}(\Gamma_{W};\k),$$
where $W=\{i\in\{1,\dots,n\}:\sigma_i\neq 0\}$.
\end{next}

We recall that a graph $G$ is a pair of sets $(V(G), E(G))$, where $E(G)$ is a familly of subsets of $V$ of cardinality $2$ called edge-set of $G$. The edge ideal of $G$ is $I(G)=(x_ix_k:\{x_i,x_k\}\in E(G))$. A {\bf clique} of $G$ is a subset $T$ of vertices of $G$ such that for all $v,w\in T$ $\{v,w\}\in E(G)$, i.e. the restriction of $G$ on the vertex subset $T$ is a complete graph. The {\bf clique complex} of $G$ is the simplicial complex $\Gamma(G)$ whose faces are the cliques of $G$. We remark that the $1$-skeleton of $\Gamma(G)$ is $G$. So, if $\Gamma$ is a simplicial complex, then $\Gamma$ is a clique complex if $\Gamma=\Gamma(\Gamma_1)$, where $\Gamma_1$ is the $1$-skeleton of $\Gamma$. We have that $I_{\Gamma}=I(G^c)$ where $G^c$ is the graph with vertex-set $V(G)$ and edge-set $E(G^c)=\{\{v,w\}\subset V(G):\{v,w\}\notin E(G)\}$, so $I_{\Gamma}$ is a square-free quadratic monomial ideal and reciprocally for any square free quadratic monomial ideal $I$ there exists a graph $G$ such that $I=I_{\Gamma(G)}$.

\igno{\begin{next}\nx{Proposition}\label{def-equiv-cliques}
	The following statements are equivalent for a simplicial complex $\Gamma$:
	\begin{enumerate}
		\item $\Gamma$ is a clique complex.
		\item Every minimal non-face of $\Gamma$ consists of $2$ vertices.
		\item $I_\Gamma$ is generated by quadratic monomials.
	\end{enumerate}
\end{next}}
\igno{\proof
\begin{itemize}
 \item[$2\Leftrightarrow 3$)] It's clear that (2) and (3) are equivalent.
 \item[$1\Rightarrow 2$)] Let be any $F\subset V(\Gamma)$, such that $F\notin\Gamma$, then the complete graph of vertices $F$ is not a subgraph of $\Gamma_1$, so there is two elemens $x\neq y \in F$ where $\{x,y\}\notin\Gamma$. By this way all minimal non-face of $\Gamma$ has two vertices.
 \item[$1\Leftarrow 2$)] Let be $F\subset V(\Gamma)$ such that $F$ is a clique of $\Gamma_1$. If $|F|\leq 2$, then $F$ is a face of $\Gamma$ by definition of $1$-skeleton. If $|F|\geq 3$, as every minimal non-face of $\Gamma$ consists of $2$ vertices, if $F$ was not a face of $\Gamma$, there would be a $F_1\subset F$, $|F_1|=2$, such that $F_1$ would not be a face of $\Gamma$, so $F$ would not be a clique of $\Gamma_1$, it is a contradiction, so $F$ is a face of $\Gamma$.
\qed
\end{itemize}}
\begin{next}\nx{Definition}\label{chemin}
{\rm A cycle $C$ of a graph $G$ is a subgraph of $G$ with vertex-set $V(C)=\{v_1,\dots,v_q\}$ and its edge-set is $E(C)=\{\{v_1,v_2\},\{v_2,v_3\},\dots,\{v_{q-1},v_{q}\},\{v_q,v_1\}\}\subset E(G)$. We call lenght of $C$ to be the cardinality of $E(C)$ and we denote it as $|C|$. We will say that $C\subset\Gamma$ is a cycle of $\Gamma$, where $\Gamma$ is a simplicial complex, if and only if $C$ is a cycle of $\Gamma_1$.}
\end{next}

We say that the cycle $C$ of lenght $>3$ of $G$ has a {\bf chord} if there is an edge $\{v_i,v_j\}\in (E(G)\setminus E(C))$. We say that the cycle $C$ is minimal\index{cycle!minimal} if $|C|>3$ and it does not have any chord, equivalently $C=G_{V(C)}$. The graph $G$ is called {\bf chordal graph}\index{graphe!de chords} if all cycle of lenght $>3$ has a chord.

\begin{next}\nx{Proposition}\label{betti-p}
 Let $\Gamma=\Gamma(G)$ be the clique complex of the graph $G$, $I=I_\Gamma$ and $p\geq 1$ a natural integer. Then:
\begin{enumerate}
 \item $p_2(I)+3$ is the minimal lenght of a minimal cycle. 
 \item Besides, if $p_{2}\neq\infty$, $\beta_{p_2(I),p_2(I)+3}(I)$ is equal to the number of minimal cycles of lenght $p_2(I)+3$ of $G$. \igno{and $\beta_{p_2(I),p_2(I)+i}(I)=0$ for all $i\geq 4$.}
\end{enumerate}
\end{next}
\proof
 \begin{enumerate}
 	\item The first statement is the theorem 2.1 of \cite{eisenchido}. 
	\item We will write $p=p_2(I)$. We define ${\Cc}=\{C\subset \Gamma: C\hbox{ }{\rm is}\hbox{ }{\rm a}\hbox{ }{\rm minimal}\hbox{ }{\rm cycle}\hbox{ }{\rm of}\hbox{ }{\rm lenght}\hbox{ }p+3\}.$
By (1) $\Cc$ is not empty, because of $p=p_2(I)$. Let $C\in{\cal C}$, due to the formula of Hochster:
$$\beta_{i,V(C)}(I)={\rm dim}\widetilde{H}_{|C|-i-2}(\Gamma_{V(C)},\k)={\rm dim}_\k \widetilde{H}_{p+3-i-2}(C,\k)={\rm dim}_\k \widetilde{H}_{p-i+1}(C,\k),$$
and $$\widetilde{H}_{p-i+1}(C,\k)=\left\{\begin{array}{ll}
						1&{\rm if}\mbox{ }p-i+1=1;\\
						0&{\rm else}.
					\end{array}\right.$$
Hence $\beta_{i,{\rm supp}(V(C))}(I)=0$ for $i\neq p$, and $\beta_{p,V(C)}(I)={\rm dim}_\k \widetilde{H}_{1}(C,\k)=1$. So
$$\beta_{p,p+3}(I)=\sum_{\begin{array}{c}
															\sigma\subset V(\Gamma)\\
															|\sigma|=p+3
													 \end{array}}{\beta_{p,\sigma}}\geq \sum_{C\in {\cal C}}\beta_{p,V(C)}=|{\cal C}|.$$
Now we have to see $\beta_{p,p+3}(I)=|{\cal C}|$. Let $\sigma\subset \{1,\dots,n\}$, $|\sigma|=p+3$. Due to the formula of Hochster $\beta_{p,\sigma}(I)=\widetilde{H}_1(\Gamma_{\sigma},\k)$. We suppose $\widetilde{H}_1(\Gamma_{\sigma},\k)\neq 0$, we will show that $\Gamma_{\sigma}$ is a minimal cycle.
\begin{itemize}
 \item $\Gamma_{\sigma}$ is a clique complex, because  $\Gamma_{\sigma}$ is a sub-complex of $\Gamma$.
 \item $\Gamma_{\sigma}$ does not have minimal cycles of lenght $\leq p+2$, by part $1$. On the other hand by \cite{spanier}, pp:139-141, $\widetilde{H}_1(\Gamma_{\sigma},\k)\neq 0$ is generated by minimal cycles of $\Gamma_\sigma$, so we conclude that $\Gamma_\sigma$ must be a minimal cycle.
\end{itemize}
\mbox{ }\\
Thus $\beta_{p,p+3}(I)\leq|{\Cc}|$, but before we had seen $\beta_{p,p+3}(I)\geq|{\Cc}|$, hence $\beta_{p,p+3}(I)=|{\Cc}|$.\igno{Par ailleurs de le corollaire \ref{minimal-deg}  $\beta_{p,p+i}(I)=0$ pour tout $i\geq 4$.}
\qed
 \end{enumerate}

The second statement of the proposition \ref{betti-p} was already proved differently by O. Fernandez and P. Gimenez in \cite{oscar}.

Also, we can see that this theorem implies the Fr\"oberg's theorem:
\begin{next}\nx{Theorem}\label{frobergt}\atrib{\cite{froberg}}
 Let $\Gamma=\Gamma(G)$ be the clique complex of the graph $G$. $I=I_\Gamma$ is $2$-linear if and only if $G$ is a chordal graph.
\end{next}

%------------------------
\subsection*{Linearly joined sequences and $2$-regularity}

Let $V=<x_1,\dots,x_n>$ be a $\k$-vector space generated by all the variables of $\si$. For all $Q\subset V$, let us denote $<Q>\subset V$ as the $\k$-vector space generated by $Q$ and $(Q)\subset\si$ the generated ideal by $Q$.\\

Let $l\in\N$ and $\J_1,\dots,\J_l$ be an ordered sequence of ideals in $\si$, such that $\J_i=(\M_i,({\cal Q}_i))$, for all $1\leq i\leq l$, where $\Qq_i\subset V$ is a $\k$-vector subspace and $\M_i$ is an ideal containing no linear forms. Assume that the intersection $\J_1\cap\dots\cap \J_l$ is not redundant. The sequence of ideals $\J_1,\dots,\J_l$ is said to be linearly joined if for all $k=2,\dots,l$:
\begin{equation}\label{def-lin-joint}
 \J_k+\bigcap_{i=1}^{k-1}\J_l=(\Qq_k)+\left(\bigcap_{i=1}^{k-1}\Qq_i\right).
\end{equation}

\begin{next}\nx{Remark}\label{joint-sequence-heritage}
 {\rm Let $\J_1,\dots,\J_l$ be a linearly joined sequence of ideals in $\si$, then the sequence $(\Qq_1),\dots,(\Qq_l)$ is also linealy joined, since:
$$(\Qq_k)+\left(\bigcap_{i=1}^{k-1}\Qq_i\right)\subset (\Qq_k)+\bigcap_{i=1}^{k-1}(\Qq_i)\subset \J_k+\bigcap_{i=1}^{k-1}\J_l=(\Qq_k)+\left(\bigcap_{i=1}^{k-1}\Qq_i\right).$$
}
\end{next}

\begin{next}\nx{Theorem}\label{2-linear-prime}\atrib{See \cite{eisenchido} and Theorem 3 \cite{chidmor}}
 Let  $\J_1,\dots,\J_l$ be a sequence of prime  homogeneous ideals. The following statements are equivalent:
	\begin {enumerate}
	 \item The ideal $\J= \J_1\cap\dots\cap\J_l$ is $2$-regular
	 \item For all $j\in\{1,\dots,l\}$, ${\rm reg}(\J_i)\leq 2$ and there exist an arrangement of the sequence $\J_1,\dots,\J_l$ such that it is linearly joined.
	\end {enumerate}
\end{next}

The next proposition follows from corollary 3 and the proposition 1 in \cite{chidmor}.

\begin{next}\nx{Proposition}\label{2-reg&joint-seq}
 Let $(\Qq_1),\dots(\Qq_l)$ be a sequence of linear ideals, such that $\Qq:=(\Qq_1)\cap\dots\cap(\Qq_l)$ is not redundant. The following statements are equivalents:
\begin{enumerate}
 \item The ideal $\Qq$ is $2$-linear.
 \item There exists a permutation of the prime components of $\Qq$, namely $(\Qq_1),\dots,(\Qq_l)$, such that this sequence is linearly joined.
 \item There exist an permutation of the prime component of $\Qq$, we say $(\Qq_1),\dots,(\Qq_l)$, and for all $i\in\{2,\dots,l\}$, there exists $\k$-linear subspaces $<\D_i>,\Pp_i$, such that:
 \begin{enumerate}
	\item $\Dd_i:=\oplus_{j=i+1}^l<\D_j>$;
	\item $\Qq_i=\Dd_i+\Pp_i$;
	\item For all $k\in\{2,\dots,l\}$, and $j<k$, $<\D_j>\times\Pp_j\subset(\Pp_k)$.
 \end{enumerate}
\end{enumerate}
\end{next}

For a more detailed exposition of linearly joined sequences we suggest to read  \cite{chidmor} and \cite{eisenchido}.
%------------------------
\subsection*{Binomial extensions of clique complexes}
Let $\Gamma$ be a clique complex with vertex set $V(\Gamma)=\{x_0,x_1,\dots ,x_n\}$ and $\FF_{\Gamma}=\{F \mbox{ }{\rm facet}\mbox{ }{\rm of}\mbox{ } \Gamma\}$.  We recall the construction done by M. Morales and studied with M.L. Ha in \cite{hlam}: let $F\subset \Gamma$ be a facet of $\Gamma$ such that $\{x_{i_0}^{(F)},x_{i_1}^{(F)}\},\{x_0^{(F)},x_{i_2}^{(F)}\},$ $\dots,\{x_0^{(F)},x_{i_{m_F}}^{(F)}\}$ are {\bf proper edges}, i.e. each edge belongs uniquely to a facet who is actually $F$. To each proper edge, for example $\{x_0^{(F)},x_j^{(F)}\}$, we associate a set $Y_j^{(F)}$, which could be $\emptyset$. We define $Y^{(F)}=\cup_{j=1}^{m_F} Y_j^{(F)}$ and $\overline{F}=F\cup Y^{(F)}$. So the extension of $\Gamma$, denoted by $\overline{\Gamma}$, is the complex whose facets are $\overline{F}$ for all $F\in \FF_{\Gamma}$. If $Y^{(F)}\neq\emptyset$ we define ${\cal I}_{F}$ as the ideal of $2\times 2$ minors of the scroll matrix:

\begin{center}
$M_F=\left(  \begin{array}{llll|lll|l|lll}
    	x_{i_0}^{(F)}&y_{i_{1} 1}^{(F)}&	\dots&	y_{i_{1} n_{i_1}}^{(F)}&y_{i_{2} 1}^{(F)}&\dots&y_{i_{2} n_{i_{2}}}^{(F)}&\dots&	y_{i_{m_F} 1}^{(F)}& \dots&	y_{i_{m_F} n_{i_{m_F}}}^{(F)}\\
	y_{i_{1} 1}^{(F)}&y_{i_{1} 2}^{(F)}&	\dots&	x_{i_{1}}^{(F)}&y_{i_{2} 2}^{(F)}&	\dots&	x_{i_{2}}^{(F)}&\dots&	y_{i_{k_F} 2}^{(F)}& \dots&	x_{i_{k_F}}^{(F)}
      \end{array}\right),$
\end{center}
otherwise $\I_F=0$. Hence we define the binomial extension $\B_{\overline{\Gamma}}\subset\k[V(\overline{\Gamma})]$ of the simplicial complex $\Gamma$ as:
\begin{center}
$\B_{\overline{\Gamma}}=\I_{\overline{\Gamma}}+\J_{\overline{\Gamma}}$, where $\J_{\overline{\Gamma}}=\sum_{F\in\FF_{\Gamma}}{\cal I}_{F}.$
\end{center}

\begin{next}\nx{Proposition}\atrib{\cite[Proposition 1.4] {morales}}\label{primary-decomposition}
 For any facet $F$ of $\Gamma$, let ${\cal J}_F=(\I_F,(V(\overline{\Gamma})\setminus\overline{F}))\subset\k[V(\overline{\Gamma})]$. The primary decomposition of $\B_{\overline{\Gamma}}$ over $\k[V(\overline{\Gamma})]$ is:
$\B_{\overline{\Gamma}}=\bigcap_{F\mbox{ }{\rm facette}\mbox{ }{\rm de}\mbox{ }\overline{\Gamma}}{\cal J}_F.$
\end{next}
%Our purpose is to find $p_2(\B_{\overline{\Gamma}})$.
%------------------------
\section{Gr\"obner basis of $\B_{\overline{\Gamma}}$}

\begin{next}\nx{Proposition}\label{ideal-deter}
\rm Let ${\bf x}=\{x_1,\dots,x_n\}$ and ${\bf y}=\{y_1,\dots,y_n\}$ be two families of variables which are not necessarily distinct, but we can order ${\bf x}\cup {\bf y}$ such that $x_i>x_j$ where $1\leq i<j\leq n$ and $x_i>y_i$ where $1\leq i\leq n $. Let $\k[{\bf x}\cup{\bf y}]$ be the polynomial ring over $\k$ with ${\bf x}\cup{\bf y}$ as indeterminates, whose monomials are ordered by the lexicographic order and $M$ be the following matrix over $\k[{\bf x}\cup{\bf y}]$

$$M=\left(  \begin{array}{llll}
    	x_{1}&x_{2}&	\dots&	x_{n}\\
	y_{{1}}&y_{2}&	\dots&	y_{n}
      \end{array}\right)$$
Let $J:=I_2(M)$ be the ideal generated for all minors of $2\times2$ over $M$. Thus
\begin{enumerate}
 \item $B=\{\underline{x_{i}y_{j}}-x_{j}y_{i}: 1\leq i<j\leq n\}$ is a minimal system of generators and $B$ is also a Gr\"obner basis of $J$, where the underlined term is the leading term.
 \item ${\rm in}(J)$ is a $2$-linear ideal, which implies $J$ is a $2$-linear ideal. 
\end{enumerate}
\end{next}
\proof
\begin{enumerate}
	\item Let $f\neq g\in B$, therefore we can write them as: $f=x_iy_j-x_{j}y_{i}$ and $g=x_ky_l-x_{l}y_{k}$, where $i<j$ and $k<l$. Hence ${\rm in}(f)=x_iy_j$ and ${\rm in}(g)=x_ky_l$. We define: $$S(f,g):=\frac{{\rm lcm}({\rm in}(f),{\rm in}(g))}{{\rm in}(f)}f-\frac{{\rm lcm}({\rm in}(f),{\rm in}(g))}{{\rm in}(g)}g,$$ 
where ${\rm lcm}({\rm in}(f),{\rm in}(g))$ is the least common multiple of ${\rm in}(f)$ and ${\rm in}(g)$. Due to the Buchberger theorem (see theorem 1.7.4 of \cite{lousta}) it is enough to prove that we can reduce $S(f,g)$ to $0$ by $B$ using the division algorithm, whenever the maximal common divisor $({\rm in}(f),{\rm in}(g))\neq 1$. Suppose $({\rm in}(f),{\rm in}(g))\neq 1$. We have four cases:
\begin{enumerate}
	\item $x_i=x_k$
	\item $x_i=y_l$
	\item [(b')] $y_j=x_k$
	\item [(c)] $y_j=y_l$
\end{enumerate}
Note that cases (b) and (b') are similars, so we will only show the cases (a),(b) and (c).
\begin{enumerate}
	\item Since $x_i=x_k$, we have that $y_i=y_k$ and $$S(f,g)=-y_lx_{j}y_{i}+y_jx_{l}y_{i}.$$ We can suppose $l<j$, let $h=x_{l}y_{j}-x_jy_{l}\in B$, so ${\rm in}(h)=x_{l}y_{j}$, and $S(f,g)=y_{i}h\in B$.
	\item Since $x_i=y_l$, $$S(f,g)=-x_kx_{j}y_{i}+y_jx_{l}y_{k},$$ and $x_k>x_l>y_l=x_i>x_j>y_j$, let $h_1=x_{l}y_{j}-x_{j}y_{l}\in B$, so ${\rm in}(h_1)=x_{l}y_{j}$, and $$S(f,g)-y_{k}h_1=y_{k}x_{j}y_{l}-x_kx_{j}y_{i}=y_{k}x_{j}x_{i}-x_kx_{j}y_{i}.$$ Moreover, $h_2=x_{k}y_{i}-y_kx_{i}\in B$ with ${\rm in}(h_2)=x_{k}y_{i}$, so $S(f,g)-y_kh_1=x_{j}h_2,$ and we have that  $$S(f,g)=y_kh_1+x_{j}h_2\in B$$
	\item Since  $y_j=y_l$, so $x_j=x_l$ and $$S(f,g)=-x_kx_{j}y_{i}+x_ix_{j}y_{k}.$$ We can suppose $i<k$, let $h=x_{i}y_{k}-x_{k}y_{i}\in B$, where ${\rm in}(h)=x_{i}y_{k}$, and $S(f,g)=x_{j}h\in B.$
\end{enumerate}
Then in all the cases we conclude that $S(f,g)$ is reduced to $0$ by $B$, for $f,g\in B$ with\break $({\rm in}(f),{\rm in}(g))\neq 1$. Hence $B$ is a Gr\"obner basis of $J$. Since all the elements of $B$ are homogeneus quadratic binomials we have that if a leading term $X$ divide a term $Y$ of a binomial of $B$, then $X=Y$. We can remark that a leading term $x_iy_j$ ($1\leq i< j\leq n$) of a binomial of $B$ is not a term of other binomial of $B$, due to the definition of $M$. Thus, $B$ is a reduce Gr\"obner basis of $J$ and since $B$ is homogeneous, $B$ is a minimal system of generators of $J$.
\item Due to (1) we have that ${\rm in}_{>{\rm lex}}(J)=(x_{i}y_{j}: 1\leq i<j\leq n)$ is a quadratic square free monomial ideal thanks of the order restrictions to the variables ${\bf x}\cup {\bf y}$. Besides, we can order the generators of  ${\rm in}(B)$ in the following diagram:
$$\begin{array}{ccccccccccccccc}
&&&&&&&x_1y_n\\

&&&&&&x_1y_{n-1}&&x_2y_n\\
&&&&&\cdot&&&&\cdot\\
&&&&\cdot&&&&&&\cdot\\
&&&\cdot&&&&&&&&\cdot\\
&&x_1y_4&&&&&&&&&&\cdot\\
&x_1y_3& &x_2y_4&&&&&&&&&&\cdot\\
 x_1y_2& &x_2y_3& &x_3y_4& &&&&&&&&&x_{n-1}y_n
\end{array}$$\\

Now, we define $\Delta_i:=\{y_i\}$ and $\Pp_i:=<x_1,\dots,x_{i-1}>$, for all $i=2,\dots,n$; 
\begin{center}
$\Dd_i:=\oplus_{j=i+1}^n<\D_j>=<y_{i+1},\dots,y_{n}>$ for all $i=1,\dots,n$ and\\ $\Qq_i:=\Dd_i\oplus\Pp_i=\{x_{1},\dots,x_{i-1},y_{i+1},\dots,y_{n}\}$, for all $i=1,\dots, n$.
\end{center}
We can see that for all $k\in\{2,\dots,l\}$ and $j<k$, $<\D_j>\times\Pp_j\subset(\Pp_k)$. Then, owing to the proposition \ref{2-reg&joint-seq} $\Qq_1,\dots,\Qq_n$ is a linearly joined sequence and $(\Qq)=\cap_{i=1}^n(\Qq_i)$ is $2$-linear. In fact, by the theorem 9 \cite{chidmor} we have that the generators of $(\Qq)$ are the elements of the triangle above, thus ${\rm in}(J)=(\Qq)$ so ${\rm in}_{>{\rm lex}}(J)$ is $2$-linear. Moreover, by the proposition \ref{betti-p} we have that 
$\beta_{i,j}(J)=\beta_{i,j}({\rm in}(J))$ for all $i,j\in\N$.
Consequently, $J$ is also $2$-linear.
\qed
\end{enumerate}

The fact that $J$ is $2$-linear follows also from  theorem 4.1 (cf. \cite{eisen-goto}) and the next proposition:
\begin{next}\nx{Proposition}\label{per-scroll}
 Let ${\bf x}=\{x_1,\dots,x_n\}$ and ${\bf y}=\{y_1,\dots,y_n\}$ be two families of variables who are not necessarily different, but we can order ${\bf x}\cup {\bf y}$ such that $x_i>x_j$ where $1\leq i<j\leq n$ and $x_i>y_i$ where $1\leq i\leq n $. We set
$$M=\left(  \begin{array}{llll}
    	x_{1}&x_{2}&	\dots&	x_{n}\\
	y_{{1}}&y_{2}&	\dots&	y_{n}
      \end{array}\right)$$
Then, by applying permutations on the columns of $M$ we obtain a scroll matrix $M'$.
\end{next}
\proof
For $k=1,\dots,n$ we will define a scroll matrix $M'_k$, such that $M_k'=[M_{k-1}',C_{k}']$.
\begin{enumerate}
 \item Let $M'_{1}=\left(\begin{array}{l}x_1\\y_1 \end{array}\right)$.
 \item Suppose that we have already built $M_k$. Let $\left(\begin{array}{l}x_{j_k}\\y_{j_k} \end{array}\right)$ the $k$-th column of $M_k'$.
	\begin{enumerate}
		\item If $y_{j_k}$ appears in the first line of $M$, then there exist $j_{k+1}\in\{1,\dots,n\}$ such that $x_{j_{k+1}}=y_{j_k}$ and we take as $k+1$-th column of $M_{k+1}'$ the column $\left(\begin{array}{l}y_{j_k}\\y_{j_{k+1}} \end{array}\right)$.
		\item If $y_{j_k}$ does not appears in the first line of $M$, then we take the first column of $M$ that is not in $M'_k$ in order to be the $k+1$-th column of $M_{k+1}'$.
	\end{enumerate}
\end{enumerate}
It is clear that $M'_{k+1}$ is a scroll matrix. We define $M'=M_n'$ and we have the assertion of the proposition.

\qed

\begin{next}\nx{Notation}\label{notation-matrice}
 {\rm Let $M$ be a matrix, we will denote for $[M]_{i,j}$ the element of $M$ that is in the $i$-row and $j$-column.}
\end{next}

\begin{next}\nx{Definition}\label{ordre-cherche}
 {\rm For $i\in\{1,\dots,k\}$ let $M_i$ be a $2\times n_i$ matrix. We say that $M_1,\dots,M_k$ are admissibly ordered if and only if for all $1\leq i\leq k$ one of the following statements is satisfied:

\begin{enumerate}
 \item  Either $[M_{F_i}]_{1,1}$ is not in the second row of $M_{F_j}$, for all $j>i$;
 \item  or there exists $j>i$ such that $[M_{F_{i}}]_{1,1}=[M_{F_{j}}]_{2,1}$ and there exists $i'<i$ such that\break $[M_{F_{i'}}]_{1,1}=[M_{F_{j}}]_{1,1}$ and for all $j'<i'$, $[M_{F_{i}}]_{1,1}\neq [M_{F_{j'}}]_{1,1}$.
\end{enumerate}
}
\end{next}

%In the next example we introduce a clique complex $\Gamma$ and a binomial extension  $\overline{\Gamma}$ of $\Gamma$, such that any order of facets of $\Gamma$ is not an admissible order: 

\begin{next}\nx{Definition}\label{allow-ad}
 {\rm Let $\Gamma$ be a clique complex and $\overline{\Gamma}$ be a binomial extension of $\Gamma$. Let $\FF=\{F \mbox{ }{\rm facet}\mbox{ }{\rm of}\mbox{ } \Gamma: F\neq\overline{F}\}$. We will say that $\FF$ is endowed with an admissible order if there is an order for $\FF$, namely, $\FF=\{F_1,\dots,F_k\}$, such that the matrices $M_{F_1},\dots,M_{F_k}$ are admissibly ordered.}
\end{next}

\begin{next}\nx{Definition}\label{permutation-permise}
 {\rm Let $M$ be a $2\times n$ matrix, as in the proposition \ref{per-scroll}. Let $\Pi\in S_n$, we define 
$$\Pi(M)=\left(  \begin{array}{llll}
    	x_{\Pi(1)}&x_{\Pi(2)}&	\dots&	x_{\Pi(n)}\\
	y_{\Pi(1)}&y_{\Pi(2)}&	\dots&	y_{\Pi(n)}
      \end{array}\right)$$
We say that the permutation $\Pi\in S_n$ is an {\bf admissible permutation} of $M$ if $\Pi(1)=1$ and for all $i\in\{2,\dots,n\}$ and $j\leq i$, $y_{\Pi(i)}\neq x_{\Pi(j)}$.
}
\end{next}

\begin{next}\nx{Definition}\label{toutes-permutation}
 {\rm Let $\Gamma$ be a clique complex, $\overline{\Gamma}$ a binomial extension of $\Gamma$ and $\FF=\{F_1,\dots,F_k\}$. Let $\Pp_{\overline{\Gamma}}$ the set: {\small $$\{({\Pi_1},\dots,\Pi_k)\in S_{|Y_{F_1}|+1}\oplus\cdots \oplus S_{|Y_{F_k}|+1}:\forall h\in\{1,\dots,k\} \mbox{ }\Pi_h\in S_{|Y_{F_h}|+1}\mbox{ } {\rm admissible}\mbox{ }{\rm permutation}\mbox{ }{\rm of}\mbox{ }M_{F_h}\}.$$}}
\end{next}

\begin{next}\nx{Definition}\label{orden1}
{\rm  Let $\Gamma(G)$ be a clique complex, $\overline{\Gamma}$ be a binomial extension of $\Gamma$ and\break $\FF=\{F \mbox{ }{\rm facet}\mbox{ }{\rm of}\mbox{ } \Gamma: F\neq\overline{F}\}$ is endowed with an admissible order, namely $\FF=\{F_1,\dots,F_k\}$. For $\Pi\in\Pp_{\overline{\Gamma}}$, let  $L_1'(\Pi_h(M_{F_h}))$ be the set of all the vertices that appear in the first row of the matrix $M_{F_h}$ and does not in the first row of a matrix of a previous facet. We define an order $>_\Pi$ on the vertex set of $\overline{\Gamma}$ in the following way: $L_1'(\Pi_h(M_{F_h}))$ is ordered by the decreassing order of the columns and if $z\in L_1'(\Pi_h(M_{F_h}))$ and $z'\in L_1'(\Pi_{h'}(M_{F_{h'}}))$ so $z>z'$ if and only if $h<h'$. The remaining vertices are ordered by an arbitrary order. In case where $\Pi_k$ is the identity in $S_{|Y_{F_h}|+1}$ for all $h\in\{1,\dots,k\}$, we will denote $>_{\Pi}$ by $>$.} %If there is not any confusion related with the order $\lambda$ we will write $>_{\Pi}$ and $>$ for $>_{\lambda_\Pi}$ and $>_\lambda$ respectively.}
\end{next}

\begin{next}\nx{Remark}\label{ordre-satisfait}
 {\rm Since $M_{F_1},\dots,M_{F_k}$ are admissibly ordered by the definition \ref{orden1}, we have in this case that for all $\Pi\in \Pp$, the order $>_{\Pi}$ of the vertices of $\overline{\Gamma}$ satisfies that for all $1\leq i\leq k$ and\break $1\leq j<s\leq |Y^{(F_i)}|+1$, $[\Pi_i(M_{F_i})]_{1,j}>_{\lambda_\Pi}[\Pi_i(M_{F_i})]_{1,s}$; and for  all $1\leq j<s\leq |Y^{(F_i)}|+1$ and $1\leq j\leq |Y^{(F_i)}|+1$, $[\Pi_i(M_{F_i})]_{1,j}>_{\lambda_\Pi}[\Pi_i(M_{F_i})]_{2,j}$.}
\end{next}

\igno{\begin{next}\nx{Example}\label{ordre-somm}
 {\rm Let $\Gamma$ be the complex generated by the facets $F_1=\{a,b,c\}$, $F_2=\{b,c,d\}$ and $F_3=\{a,e,f\}$, who is evidently a clique complex with generator graph $G$ like in the figure \ref{order-exem}. In the same figure we see that its proper edges are: $\{a,b\}$, $\{a,c\}$, $\{b,d\}$, $\{c,d\}$, $\{a,e\}$, $\{a,f\}$, $\{e,f\}$. Let $M_{F_1}$, $M_{F_2}$ and $M_{F_3}$ be the following matrices:\\

\begin{center}
$M_{F_1}=\left(  \begin{array}{lll|l}
    	a&u&v&w\\
	u&v&b&c
      \end{array}\right)$,
$M_{F_2}=\left(  \begin{array}{ll|l}
    	d&x&y\\
	x&b&c
      \end{array}\right)$ and
$M_{F_3}=\left(  \begin{array}{ll}
    	a&z\\
	z&f
      \end{array}\right)$
\end{center}

\begin{figure}[htb]  
  \centering
  \scalebox{1}{
  \includegraphics[viewport=200 580 380 650]{ejemorden.pdf}}
  \vspace{3mm}  
  \caption{} 
  \label{order-exem}
\end{figure}

Thus $\overline{\Gamma}$ is generated by the facets $\overline{F}_1=\{a,b,c,w,x,y\}$, $\overline{F}_2=\{b,c,d,z,y\}$ and $F_3\{a,e,f,z\}$ where the generator graph is $\overline{G}$, as we can see in the figure \ref{order-exem}. So $>$  is an order to the vertices of $\overline{\Gamma}$, and it is: $a>u>v>w>d>x>y>z>b>c>f>e$.}\qeds
\end{next}}

\begin{next}\nx{Definition}\label{complexe-sans-diagonal}
 {\rm Let $\Gamma$ be a clique complex, $\overline{\Gamma}$ a binomial extension of $\Gamma$. Let $\Pi\in\Pp_{\overline{\Gamma}}$. We define the simplicial complex $\overline{\Gamma}^{(\Pi)}$  as the clique complex generated by the graph $(\overline{\Gamma}^{(\Pi)})_1$ obtained from $(\overline{\Gamma})_1$ deleting all the edges which correspond to all the diagonals of $\Pi(M_{F})$ from the top to the bottom from left to right, where $F$ is a facet of $\Gamma$ such that $F\neq \overline{F}$. If $\Pi$ is the identity we will write $\overline{\Gamma}'$ instead of $\overline{\Gamma}^{(\Pi)}$. }
\end{next}

\begin{next}\nx{Remark}\label{restiction}
 {\rm Thanks to the definition \ref{complexe-sans-diagonal}, for any facet $F\in\Gamma$, the restriction $\overline{\Gamma}^{(\Pi)}_{\overline{F}}$ coincide with $<F>$ if $F=\overline{F}$, and coincide with $\Gamma({\rm in}(\I_F))$ if $F\neq\overline{F}$. As $\I_F$ for all $F\neq\overline{F}$ is a scroll ideal, $\I_F$ is $2$-linear and by Fr\"oberg's (theorem \ref{frobergt}) we have that $\Gamma({\rm in}(\I_F))_1$ is a chordal graph. If $F=\overline F$, $\overline{\Gamma}^{(\Pi)}_{\overline{F}}=<F>$ is the simplex $F$ so $\overline{\Gamma}^{(\Pi)}_{\overline{F}}$ is a chordal graph. Thus, for all $F$ facet of $\Gamma$, $(\overline{\Gamma}^{(\Pi)}_{\overline{F}})_1$ is a chordal graph.}
\end{next}

Let $\Gamma$ be a clique complex; $\overline{\Gamma}$ a binomial extension of $\Gamma$. Set:
${\rm NF}_{\overline{\Gamma}}$ the generator set of $\I_{\overline{\Gamma}}$, i.e.
\begin{center}
${\rm NF}_{\overline{\Gamma}}=\{{\bf x}^{{\rm supp}(\sigma)}| \sigma\notin \overline{\Gamma}\hbox{ }{\rm and}\hbox{ }|\sigma|=2\}$\\
\end{center}
and ${\rm J}_{F}$ the set of $2\times2$ minors  of the matrix $M_F$. Then
$B_{\overline{\Gamma}}={\rm NF}_{\overline{\Gamma}}\bigcup(\cup_{k=1}^{l}{\rm J}_{F_k})$ is a generator system of $B_{\overline{\Gamma}}$.

\begin{next}\nx{Proposition}\label{preliminar-in}
Let $\Gamma$ be a clique complex, $\overline{\Gamma}$ a binomial extension of $\Gamma$. We suppose\break $\FF=\{F \mbox{ }{\rm facet}\mbox{ }{\rm of}\mbox{ } \Gamma: F\neq\overline{F}\}$ is endowed with an admissible order. Let $\Pi\in\Pp_{\overline{\Gamma}}$. We order the vertices of $\overline{\Gamma}$ by $>_{\Pi}$ (see definition \ref{orden1}). Then \begin{enumerate}
\item  $B_{\overline{\Gamma}}$ is a Gr\"obner basis of $\B_{\overline{\Gamma}}$ to the lexicographic order on the monomials of $\sis$.
\item  The ideal ${\rm in}_{>_{\Pi{\rm lex}}}(\B_{\overline{\Gamma}})$ is a square free quadratic monomial ideal and the associated simplicial complex of ${\rm in}_{>_{\Pi{\rm lex}}}(\B_{\overline{\Gamma}})$ is $\overline{\Gamma}^{(\Pi)}$.
\end{enumerate}
\end{next}
\proof
\begin{enumerate} 
\item If ${\cal I}_F$ is not trivial, then ${\cal I}_F$ is generated by the minors of $\Pi(M_F)$, i.e., $$[\Pi(M_F)]_{1,i}[\Pi(M_F)]_{2,k}-[\Pi(M_F)]_{1,k}[\Pi(M_F)]_{2,i}$$ with $i<k$, where $[\Pi(M_F)]_{1,i}[\Pi(M_F)]_{2,k}$ is the initial term for the order $\succ_\Pi$.\\
By using the Buchberger's criteria we need to prove that $$S(f,g)=\frac{l.c.m({\rm in}(f),{\rm in}(g))}{{\rm in}(f)}f-\frac{l.c.m({\rm in}(f),{\rm in}(g))}{{\rm in}(g)}g\xrightarrow{B_{\overline{\Gamma}}} 0$$ for all $f,g\in B_{\overline{\Gamma}}$ whenever $({\rm lm}(f),{\rm lm}(g))\neq 1$. We have four cases:
\begin{enumerate}
	\item $f,g\in{\rm NF}_{\overline{\Gamma}}$, so they are monomials and $S(f,g)=0$
	\item $f,g\in {\rm J}_{F}$ for any facet $F$ of $\Gamma$. As $f,g$ are minors of the matrix $M_F$
and $[\Pi(M_F)]_{1,i}>_{\lambda_\Pi}[\Pi(M_F)]_{1,j}$, where $1\leq i<j\leq |Y^{(F)}|+1$ and $[\Pi(M)]_{1,i}\succ[\Pi(M)]_{2,i}$ for\break$i\in\{1,\dots,|Y^{(F)}|+1\}$; thus, by the proposition \ref{ideal-deter} ${\rm J}_{F}$ is a Gr\"obner basis of $\I_F$, and then $S(f,g)\xrightarrow{J_{F}} 0$.
	\item $f\in{\rm NF}_{\bar{\Gamma}}$ and $g\in {\rm J}_{F}$ for some facet $F$ of ${\Gamma}$. Let $$g=[\Pi(M_F)]_{1,i}[\Pi(M_F)]_{2,k}-[\Pi(M_F)]_{1,k}[\Pi(M_F)]_{2,i},$$ where $i<k$, and $f=ab$. Hence, we can suppose that $a=[\Pi(M_F)]_{1,i}$ or $a=[\Pi(M_F)]_{2,k}$. So $b\notin F$ and $S(f,g)=b[\Pi(M_F)]_{1,k}[\Pi(M_F)]_{2,i}.$ Besides $[\Pi(M_F)]_{1,k}\in Y^{(F)}$. Consequently $b[\Pi(M_F)]_{1,k}\in {\rm NF}_{\overline{\Gamma}}$ and $S(f,g)\rightarrow 0$.
	\item $f\in{\rm J}_{F}$ and $g\in{\rm J}_{F'}$ for two different facets $F$ and $F'$ of $\Gamma$. 
Since $({\rm in}(f),{\rm in}(g))\neq 1$, the matrices $M_{F}$ and $M_{F'}$ must have a common vertex $x$ with $x\in F\cap F'$ such that $x|{\rm in}(f)$ and $x|{\rm in}(g)$. We have four possibilities for this vertex $x$:
\begin{enumerate}
	\item Either $x=[\Pi(M_F)]_{1,1}=[\Pi(M_{F'})]_{1,1}$.\\
	Hence $f=x[\Pi(M_F)]_{2,k}-[\Pi(M_F)]_{1,k}[\Pi(M_F)]_{2,1}$ and $g=x[\Pi(M_{F'})]_{2,r}-[\Pi(M_{F'})]_{1,r}[\Pi(M_{F'})]_{2,1}$ for some $k\in\{1,\dots,|Y^{(F)}|+1\}$  and some $r\in\{1,\dots,|Y^{(F')}|+1\}$. Thus $$S(f,g)=-[\Pi(M_{F'})]_{2,r}[\Pi(M_F)]_{1,k}[\Pi(M_F)]_{2,1}+[\Pi(M_{F})]_{2,k}[\Pi(M_{F'})]_{1,r}[\Pi(M_{F'})]_{2,1}.$$
	For $[\Pi(M_F)]_{2,k}$ we have to cases:
	\begin{enumerate}
	 \item $[\Pi(M_F)]_{2,k}\in V(\Gamma)$. Then $[\Pi(M_F)]_{2,k}\notin\overline{F'}$. Otherwise $\{x,[\Pi(M_F)]_{2,k}\}\in \overline{F'}$ and $\{x,[\Pi(M_F)]_{2,k}\}$ would not be a proper edge of $F$.
	\item $[\Pi(M_F)]_{2,k}\notin V(\Gamma)$. Then, by definition of $\B_{\overline{\Gamma}}$, $[\Pi(M_F)]_{2,k}\notin\overline{F'}$.
	\end{enumerate} 
 	So, in any case $[\Pi(M_F)]_{2,k}\notin\overline{F'}$, and by the same way we can see that $[\Pi(M_{F'})]_{2,r}\notin \overline{F}$. Moreover, as $k,r\geq 2$ we have that $[\Pi(M_F)]_{1,k},[\Pi(M_{F'})]_{1,r}\notin V(\Gamma)$. Consequently $$[\Pi(M_{F'})]_{1,r}[\Pi(M_F)]_{1,k},\mbox{ }[\Pi(M_{F})]_{2,k}[\Pi(M_{F'})]_{1,t}\in {\rm NF}_{\bar{\Delta}},$$ hence:
$S(f,g)\xrightarrow{\{[\Pi(M_{F'})]_{2,r}y_{i_1}^{(F)},[\Pi(M_{F})]_{2,k}y_{j_1 1}^{(F')}\in {\rm NF}_{\bar{\Delta}}\}} 0.$
	\item Either $x=[\Pi(M_F)]_{1,1}=[\Pi(M_{F'})]_{2,t}$, where $t\in\{1,\dots,|Y^{(F')}|+1\}$.\\
		Moreover $f=x[\Pi(M_F)]_{2,k}-[\Pi(M_F)]_{1,k}[\Pi(M_F)]_{2,1}$ for some $k\in\{2,\dots,|Y^{(F)}|+1\}$ and\break $g=[\Pi(M_{F'})]_{1,r}x-[\Pi(M_{F'})]_{1,t}[\Pi(M_{F'})]_{2,r}$ for some $r\in\{1,\dots,|Y^{(F')}|+1\}$.
		Since\break ${\rm in}(g)=[\Pi(M_{F'})]_{1,r}x$, we have that  $t\geq 2$ and
		$$S(f,g)=-[\Pi(M_{F'})]_{1,r}[\Pi(M_F)]_{1,k}[\Pi(M_F)]_{2,1}^{(F)}+[\Pi(M_{F})]_{2,k}[\Pi(M_{F'})]_{1,t}[\Pi(M_{F'})]_{2,r}.$$
		As in the case i. we can prove that $[\Pi(M_F)]_{2,k}\notin\overline{F'}$ and $[\Pi(M_{F'})]_{1,r}\notin \overline{F}$.\\
		Moreover, as $k,t\geq 2$ we have that $[\Pi(M_F)]_{1,k},[\Pi(M_{F'})]_{1,t}\notin V(\Gamma)$. Therefore
		$$S(f,g)\xrightarrow{[\Pi(M_{F'})]_{2,r}[\Pi(M_F)]_{1,k},[\Pi(M_{F})]_{2,k}y_{r,n_r}^{(F')}}0.$$
	%\item $x=x_{i_p}^{(F)}=x_{i_0}^{(F')}$, o\`u $p\in\{1,\dots,m_F\}$ est similaire au cas $b.$
	\item Or $x=[\Pi(M_F)]_{1,t}=[\Pi(M_{F'})]_{2,t'}$, where $t\in\{2,\dots,|Y^{(F)}|+1\}$ and\break$t'\in\{2,\dots,|Y^{(F')}|+1\}$.\\
		Thus $f=[\Pi(M_F)]_{1,k}x-[\Pi(M_{F})]_{1,t}[\Pi(M_F)]_{2,k}$ for some $k\in\{1,\dots,|Y^{(F)}|+1\}$ and
		\begin{center}
		$g=[\Pi(M_{F'})]_{1,r}x-[\Pi(M_{F'})]_{1,t'}[\Pi(M_{F'})]_{2,r}$ for some $r\in\{1,\dots,|Y^{(F')}|+1\}$.		\end{center}
		 As in the last case we can prove that $t\geq 2$ and $t'\geq 2$. So
		$$S(f,g)=-[\Pi(M_{F'})]_{1,r}[\Pi((M_{F})]_{1,t}[\Pi(M_F)]_{2,k}+[\Pi(M_{F})]_{1,k}[\Pi(M_{F'})]_{1,t'}[\Pi(M_{F'})]_{2,r}.$$ 
		And like in the other cases we can prove that $[\Pi(M_F)]_{1,k}\notin\overline{F'}$ and $[\Pi(M_{F'})]_{1,r}\notin \overline{F}$.\\
		Besides, since $t,t'\geq 2$ we get that $[\Pi(M_F)]_{1,t},[\Pi(M_{F'})]_{1,t'}\notin V(\Gamma)$. 
		Then $$[\Pi(M_{F'})]_{1,r}[\Pi(M_{F})]_{1,t},\mbox{ }[\Pi(M_F)]_{1,k}x[\Pi(M_{F'})]_{1,t'}\in {\rm NF}_{\bar{\Delta}},$$
		 and $S(f,g)\xrightarrow{[M_{F'}]_{1,r}[M_{F}]_{1,t},[M_F]_{1,k}x[M_{F'}]_{1,t'}}0.$
\end{enumerate}
\end{enumerate}
Thus for all $f,g\in\B_{\overline{\Gamma}}$, $S(f,g)\xrightarrow{\B_{\bar{\Gamma}}}_+ 0$, and by the Buchberger algorithm we can conclude $B_{\overline{\Gamma}}^{(\Pi)}$ is a Gr\"obner basis of $\B_{\overline{\Gamma}}$.
\item As for all matrix $M_F$, $${\rm in}([\Pi(M_F)]_{1,i}[\Pi(M_F)]_{2,k}-[\Pi(M_F)]_{1,k}[\Pi(M_F)]_{2,i})=[\Pi(M_F)]_{1,i}[\Pi(M_F)]_{2,k}$$ with $i<k$ and $[\Pi(M_F)]_{1,i}\neq[\Pi(M_F)]_{2,k}$, we have that ${\rm in}(B)$ is a square-free monomial set. And by $1.$ $B$ is a Gr\"obner basis of $\B_{\Gamma}$, so $\B_{\overline{\Gamma}}={\rm in}(B)$ is a square-free monomial ideal and by definition of $\Gamma({\rm in}(\B_{\overline{\Gamma}}))$ we get that $\Gamma({\rm in}(\B_{\overline{\Gamma}}))=\overline{\Gamma}^{(\Pi)}$.
\qed
\end{enumerate}
%------------------------
\section{Existence of admissible orders}

It's not obvious that given a graph $G$, $\Gamma(G)$ its clique complex and $\overline{\Gamma}$ a simplicial extension of $\Gamma$,\break $\FF=\{F \mbox{ }{\rm facet}\mbox{ }{\rm of}\mbox{ } \Gamma: F\neq\overline{F}\}$ is endowed with an admissible order; as we can see in the next example:

\begin{next} \nx{Example}\label{extensiones}
{\rm Let $\Gamma$ be the clique complex generated by the graph $G$ of the figure \ref{malorder-exem}. Then,\break $F_1=\{a,e,b\}$, $F_2=\{d,h,a\}$, $F_3=\{c,g,d\}$, $F_4=\{b,f,c\}$ are the facets of $\Gamma$. We associate to each facet the following matrices:

\begin{center}
$M_{F_1}=\left(  \begin{array}{lll|l}
    	a&x&y&z\\
	x&y&b&e
      \end{array}\right)$,
$M_{F_2}=\left(  \begin{array}{lll|l}
    	d&r&s&q\\
	r&s&a&h
      \end{array}\right)$,

\end{center}
\begin{center}
$M_{F_3}=\left(  \begin{array}{ll|l}
    	c&u&t\\
	u&g&d
      \end{array}\right)$ and
$M_{F_4}=\left(  \begin{array}{ll|l}
    	b&w&v\\
	w&c&f
      \end{array}\right)$,
\end{center}

\begin{figure}[htb]  
  \centering
  \scalebox{0.70}{
  \includegraphics[viewport=150 520 320 630]{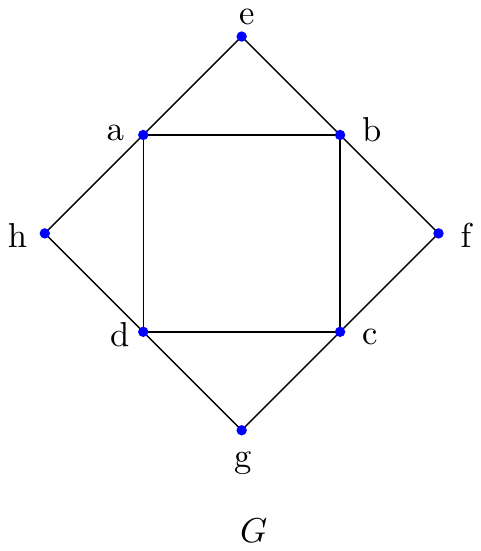}}
  \vspace{3mm}  
  \caption{} 
  \label{malorder-exem}
\end{figure}

Note that for any possible order $\lambda$ of the matrices $M_{F_1},M_{F_2},M_{F_3},M_{F_4}$, namely $\lambda=M_1,M_2,M_3,M_4$, thus there exists $i\in\{2,3,4\}$ such that $[M_1]_{1,1}$ appears in the second row of $M_i$. Nevertheless if we replace $M_{F{2}}$ by 
$$M'_{F_2}=\left(  \begin{array}{lll|l}
    	a&s&r&h\\
	s&r&d&q
      \end{array}\right)$$
we have that $M_{F_1},M_{F_4},M_{F_3},M'_{F_2}$ are admissibly ordered and we can consider that $\B_{\overline{\Gamma}}$ is the binomial extension of $I_{\Gamma'}$, where $\Gamma'$ is the clique complexe generated by the graph obtained from $G$ replacing $q$ to $h$.}
\end{next}

\begin{next}\nx{Remark}\label{ciclillo}
 {\rm Let $C_n$ be the cycle of $n$-vertices $x_1,\dots,x_n$ and $\Gamma$ the clique complex generated by $C_n$. As $$\FF=\{F \mbox{ }{\rm facet}\mbox{ }{\rm of}\mbox{ } \Gamma: F\neq\overline{F}\}=\{\{x_1,x_2\},\{x_2,x_3\},\dots,\{x_{n},x_1\}\},$$ any facet of $C_n$ is itself a proper edge. In order to obtain a binomial extension of $C_n$, for any facet $e_i=\{x_i,x_{i+1}\}$ of $C_n$, for all $i\in\{1,\dots,n\}$, where $x_{n+1}=x_1$, we can choose as $x_0^{i}$ either $x_i$ or $x_{i+1}$. So we can suppose that for all $i=1,\dots,n-1$: $M_{e_i}=0$ if $Y^{(e_i)}=0$; 
$$M_{e_i}=\left(  \begin{array}{llll}
    	x_{i}&y_{1}^{(i)}&	\dots&	y_{s_{i}}^{(i)}\\
	y_{1}^{(i)}&y_{2}^{(i)}&	\dots&	x_{i+1}
\end{array}\right),$$
if $Y^{(e_i)}\neq \emptyset$. If $s_n=0$ then $M_{e_n}=0$; and if $Y^{(e_n)}=0$, then
$$M_{e_n}=\left(  \begin{array}{llll}
    	x_{1}&y_{1}^{(n)}&	\dots&	y_{s_{n}}^{(n)}\\
	y_{1}^{(n)}&y_{2}^{(n)}&	\dots&	x_{n}
\end{array}\right).$$
So the family ${\cal M}=\{M_{e_{j_i}}\neq 0: j_i\in\{1,\dots,n\}\}$ is admissibly ordered. Finally we conclude $\FF$ is endowed by an admissible order.}
\end{next}

In the next lines we will see that there exists another family of clique complexes $\Gamma$ distinct to $\{C_n:n\in\N\}$ which satifies that $\FF$ is endowed by an admissible order.

\begin{next}\nx{Remark}\label{clique-ord-admis-trivial}
{\rm Clearly, if $|\FF|=1$, $\FF$ is endowed with an admissible order. In addition, if $\FF$ is endowed by an admissible order, then any subfamily $\FF'\subset\FF$ is endowed with an admissible order, because the admissible order of $\FF$ restricted over the subfamily $\FF'\subset\FF$ is an admissible order of $\FF'$.}
\end{next}

\begin{next}\nx{Notation}
 Let $F\neq \overline{F}$ facet of $\Gamma$, we set $L_2(M_{F_i})$ the set of all the elements of the second line of $M_{F_i}$.
\end{next}

\begin{next}\nx{Lemma}\label{cliques-ord-admiss}
Let $\Gamma$ be a clique complex and $\overline{\Gamma}$ a binomial extension of $\Gamma$. Let  $\FF'$ be an ordered subfamily, that is $\FF'=\{F_1,\dots, F_s\}\subset F$ where $s>1$, such that $[M_{F_i}]_{1,1}\in L_2(M_{F_{i+1}})$ for all $i\in \{1,\dots,s-1\}$ and $[M_{F_s}]_{1,1}\in L_2(M_{F_1})$, then $\FF '$ is not endowed with any admissible order.
\end{next}
\proof
Let $\Pi\in S_s$ any permutation of $\{1,\dots,s\}$, then $M_{F_{\Pi(1)}},\dots,M_{F_{\Pi(s)}}$ are not admissibly ordered, because:
\begin{itemize}
	\item If $\Pi(1)\neq s$, then $[M_{F_{(\Pi(1))}}]_{1,1}\in L_2(M_{F_{\Pi(1)+1}})$, thus there exists $1<j\leq s$ such that\break$\Pi(j)=\Pi(1)+1$; and $M_{F_{\Pi(1)}},\dots,M_{F_{\Pi(s)}}$ are not admissibly ordered.
	\item If $\Pi(1)=s$, then there exists $1<j\leq s$ such that $\Pi(j)=1$ and $$[M_{F_{(\Pi(1))}}]_{1,1}\in L_2(M_{F_{1}})=L_2(M_{F_{\Pi(j)}});$$ $F_{\Pi(1)},\dots,F_{\Pi(s)}$ are not admissibly ordered.
\end{itemize}
In this manner we conclude that $\FF'$ is not endowed with an admissible order.
\qed

\begin{next}\nx{Proposition}\label{clique-ord-admiss2}
Let $\Gamma$ be a clique complex generated by $G$ and $\overline{\Gamma}$ a binomial extension of $\Gamma$. $\FF$ is not endowed with an admissible order if and only if there exists a subfamily $\FF '$ such that $\FF'$ has an order $\FF'=\{F_1,\dots, F_s\}$ where $[M_{F_i}]_{1,1}\in L_2(M_{F_{i+1}})$ for all $i\in \{1,\dots,s-1\}$ and $[M_{F_s}]_{1,1}\in L_2(M_{F_1})$.
\end{next}
\proof
\begin{itemize}
 \item [$\Rightarrow)$] Let $1\leq s\leq k$ be the integer satisfying that for any subfamily $\FF^*\subset \FF$ with $|\FF^*|<s$, $\FF^*$ is endowed with an admissible order and there exists a subfamily $\FF '\subset \FF$ with $|\FF '|=s$ is not endowed with an admissible order. Obviously $s>1$. We are going to prove that for the family $\FF '$ there exists an order $\FF '=\{F_1,\dots, F_s\}$ such that $[M_{F_i}]_{1,1}\in L_2(M_{F_{i+1}})$ for all $i\in \{1,\dots,s-1\}$ and $[M_{F_s}]_{1,1}\in L_2(M_{F_1})$.
 
	\begin{enumerate}
	 \item Assume that there exists $M_F\in\FF'$ such that for all $M_{F'}\in\FF '$, $F'\neq F$,  $$[M_{F}]_{1,1}\notin L_{2}(M_{F'}).$$ By hypothesis, $\FF '\setminus\{F\}$ is endowed with an admissible order, namely $\FF'\setminus\{F\}=\{F_{1},\dots,F_{s-1}\}$. Then $[M_{F}]_{1,1}\notin L_{2}(M_{F_i})$ for all $i=1\dots, s-1$, and in this way $M_F,M_{F_{i_1}},\dots,\dots,M_{F_{i_{s-1}}}$ are admissibly ordered, but it is a contradiction because we had supposed that $\FF '$ was not endowed with an admissible order. Consequently: $\forall F\in\FF'$ there exists ${F'}\in\FF '$, $F'\neq F$ such that $[M_{F}]_{1,1}\in L_{2}(M_{F'})$.
	\item Now we are going to prove that for the family $\FF '$ there exists an order $\FF '=\{F_1,\dots, F_s\}$ such that $[M_{F_i}]_{1,1}\in L_2(M_{F_{i+1}})$ for all $i\in \{1,\dots,s-1\}$ and $[M_{F_s}]_{1,1}\in L_2(M_{F_1})$.\\
	Let any $F_1\in\FF'$, by the point $1$, there exists $F_2\neq F_1$ such that $[M_{F_1}]_{1,1}\in  L_{2}(M_{F_{2}})$. Moreover by the same point $1$, there exists $F\in(\FF'\setminus F_2)$ such that  $[M_{F_2}]_{1,1}\in  L_{2}(M_{F})$.
		\begin{itemize}
			\item [\dtb] If $F=F_1$, then the family of facets $\{F_{1},F_{2}\}$ is not endowed with an admissible order, and by the lemma \ref{cliques-ord-admiss} and by maximality of $s$, $s=2$ and $\FF'=\{F_{1},F_{2}\}$ satisfies that $[M_{F_i}]_{1,1}\in L_2(M_{F_{2}})$ and $[M_{F_2}]_{1,1}\in L_2(M_{F_1})$.
			\item [\dtb] If $F\neq F_1$, $s>2$, we call $F_3:=F$. Thanks to $1.$ there exists $F\in\FF'\setminus F_3$ such that  $[M_{F_2}]_{1,1}\in  L_{2}(M_{F_{3}})$.
				\begin{itemize}
					\item If $F=F_1$, then the family of facets $\{{F_{1}},F_{2},F_{3}\}$ is not endowed with an admissible order, and by the lemma \ref{cliques-ord-admiss} and by maximality of $s$, $s=3$ and the family $\FF'=\{F_{1},F_{2},F_{3}\}$ is endowed by an admissible order that we wanted to prove.
					\item If $F=F_2$, the subfamily of facets $\{F_{2},F_{3}\}$ does not allow and admissible order by the lemma \ref{cliques-ord-admiss} but it is not possible because $s>2$ is maximal.
					\item Thus if $F\neq F_1$, then $F\neq F_2$ and $s>3$. We set $F_4=F$ and by induction we can suppose that we have found a sequence of differents facets in $\FF'$: $F_1,F_2,F_3,\dots,F_{s}$ such that $[M_{F_i}]_{1,1}\in  L_{2}(M_{F_{i+1}})$ for all $i=1,\dots,s-1$. Thanks to item $1.$ there exists $F\in\FF'\setminus F_{s-1}$ such that  $[M_{F_1}]_{1,1}\in  L_{2}(M_{F})$. If $F\neq F_1$, we can say $F=F_j$ with $j>2$, then $F_j,F_{j+1},\dots, F_{s}$ satisfies that  $[M_{F_i}]_{1,1}\in L_2(M_{F_{i+1}})$ for all $i\in \{j,\dots,s-1\}$ and $[M_{F_s}]_{1,1}\in L_2(M_{F_j})$. By the lemma \ref{cliques-ord-admiss} the family $\FF''=\{F_j,F_{j+1},\dots, F_s\}$ is not endowed with an admissible order and it is a contradiction to the maximality of $s$. Then $F=F_1$, and $\FF'$ satisfies that there exists an order $\FF '=\{F_1,\dots, F_s\}$ such that $[M_{F_i}]_{1,1}\in L_2(M_{F_{i+1}})$ for all $i\in \{1,\dots,s-1\}$ and $[M_{F_s}]_{1,1}\in L_2(M_{F_1})$.
				\end{itemize}
		\end{itemize}

	\end{enumerate}
 \item [$\Leftarrow )$] By the lemma \ref{cliques-ord-admiss} $\FF '$ is not endowed with an admissible order. Thus, by the remark \ref{clique-ord-admis-trivial}, $\FF$ is not endowed with an admissible order.
\qed\\
\end{itemize}
\begin{next}\nx{Lemma}\label{cliques-ord-admis-petit}
Let $G$ be a graph, $\Gamma:=\Gamma(G)$ the clique complex generated by $G$ and $\overline{\Gamma}$ a binomial extension of $\Gamma$. If $|\FF|\leq 3$, then $\FF$ is endowed with an admissible order.
\end{next}
\proof
\begin{itemize}
 \item ${\bf k=2}$. If the family $\FF$ doesn't allow any admissible order, by the proposition \ref{clique-ord-admiss2} there exist an order of $\FF$: $\FF=\{F_1,F_2\}$ such that $[M_{F_1}]_{1,1}\in L_2(M_{F_2})$ and $[M_{F_2}]_{1,1}\in L_2(M_{F_1})$. Thus $\{[M_{F_1}]_{1,1},[M_{F_2}]_{1,1}\}\in F_1\cap F_2$ and $\{[M_{F_1}]_{1,1},[M_{F_2}]_{1,1}\}$ is not a proper edge of $\Gamma$, but it is not possible because $\overline{\Gamma}$ is a binomial extension of $\Gamma$. Consequently $\FF$ is endowed with an admissible order.
 \item ${\bf k=3}$. If the family $\FF$ is not endowed with an admissible order, then by the proposition \ref{clique-ord-admiss2} there exist an order of $\FF$: $\FF=\{F_1,F_2,F_3\}$ such that $[M_{F_1}]_{1,1}\in L_2(M_{F_2})$, $[M_{F_2}]_{1,1}\in L_2(M_{F_3})$ and\break$[M_{F_3}]_{1,1}\in L_2(M_{F_1})$. So we can conclude that $e_1=\{[M_{F_2}]_{1,1},[M_{F_1}]_{1,1}\}\in F_2$,\break $e_2=\{[M_{F_3}]_{1,1},[M_{F_2}]_{1,1}\}\in F_3$, $e_3=\{[M_{F_1}]_{1,1},[M_{F_3}]_{1,1}\}\in F_1$ and\break$L=\{[M_{F_1}]_{1,1},[M_{F_2}]_{1,1},[M_{F_3}]_{1,1}\}$ is a clique of $G$, thus $L\subset F_i$ for some $i\in\{1,2,3\}$ and in this way either $e_1$, either $e_2$ or $e_3$ is not a proper edge, and it is a contradiction to the fact that $\overline{\Gamma}$ is a binomial extension of $\Gamma$. Thus, we conclude that $\FF$ is endowed with an admissible order.
\qed
\end{itemize}
\mbox{ }\\
\begin{next}\nx{Proposition}\label{cliques-prop-admi}
Let $G$ be a graph, $\Gamma:=\Gamma(G)$ the clique complex generated by $G$ and $\overline{\Gamma}$ a binomial extension of $\Gamma$. Assuming that the graph $G$ restricted on the family $\FF$ is a chordal graph, then $\FF$ is endowed with an admissible order.
\end{next}
\proof
If $|\FF|=1,2$ or $3$ by the lemma \ref{cliques-ord-admis-petit} the family $\FF$ is endowed by an admissible order.\\

Let us suppose that $|\FF|\geq 4$. Assume that the family $\FF$ is not endowed with an admissible order. So, by the proposition \ref{clique-ord-admiss2} and lemma \ref{cliques-ord-admis-petit} there exists a subfamily $\FF '=\{F_1,\dots, F_s\}\subset \FF$ where $s>3$, such that $[M_{F_i}]_{1,1}\in L_2(M_{F_{i+1}})$ for all $i\in \{1,\dots,s-1\}$ and $[M_{F_s}]_{1,1}\in L_2(M_{F_1})$. Thus, by definition of binomial extension, $\{[M_{F_{i+1}}]_{1,1},[M_{F_{i}}]_{1,1}\}$ is a proper edge of $F_{i+1}$ for all $i\in\{1,\dots,s-1\}$ and $\{[M_{F_{s}}]_{1,1},[M_{F_{1}}]_{1,1}\}$ is a proper edge of $F_{1}$. So, $\{\{[M_{F_{1}}]_{1,1},[M_{F_{2}}]_{1,1}\},\dots,\{[M_{F_{s}}]_{1,1},[M_{F_{1}}]_{1,1}\}\}$ are the edges of a cycle $C$ of $G$, but since $G$ restricted on the family $\FF$ is a chordal graph, $C$ is not a minimal cycle and $C$ must have a clique $V\subset V(C)$ of size $3$ with two edges $e,e'$ of $C$. In this way we have that $V\subset F_i$ for some $i\in\{1,\dots, k\}$ and at least an edge among $e$ and $e'$ is in two facets of $\Gamma$, but this is a contradiction to the fact that $e,e'$ are proper edges of $\Gamma$. Consequently $\FF$ is endowed with an admissible order.
\qed

%------------------------
\section{Lower bounds of $p_2(\B_{\overline{\Gamma}})$}\label{bornes-inf}

Let $\Gamma$ be a clique complex and $\overline{\Gamma}$ be a binomial extension of $\Gamma$. We suppose that $\FF=\{F \mbox{ }{\rm facet}\mbox{ }{\rm of}\mbox{ } \Gamma:\overline{F}\neq F\}$ is endowed with an admissible order, namely $\FF=\{F_1,\dots,F_k\}$. Let $\Pi\in \Pp_{\overline{\Gamma}}$ (see definition \ref{toutes-permutation}). Using the definition \ref{orden1} we get and order $>_\Pi$ for the vertex-set of $\overline{\Gamma}$ and by the proposition \ref{preliminar-in} the associated simplicial complex of  ${\rm in}_{>_{\Pi}{\rm lex}}(B_{\overline{\Gamma}})$ is $\overline{\Gamma}^{(\Pi)}$ defined in \ref{complexe-sans-diagonal}. Using the proposition \ref{igualdad-ini}  we have that
\begin{equation}\label{max-eq}
 {\rm max}_{\Pi\in\Pp_{\overline{\Gamma}}}(p_2(I(\overline{\Gamma}^{(\Pi)})))\leq p_2(B_{\overline{\Gamma}}).
\end{equation}
And from proposition \ref{betti-p}, let $\Cc_{\Pi}'$ be the family of all minimal cycles of $\overline{\Gamma}^{\Pi}$ 
\begin{equation}\label{max-eq2}
 {\rm max}_{\Pi\in\Pp_{\overline{\Gamma}}}(p_2(I(\overline{\Gamma}^{(\Pi)})))+3= {\rm max}_{\Pi\in\Pp_{\overline{\Gamma}}}(min_{C'\in\Cc_{\Pi}'}|C'|)
\end{equation}

So we would like to know the family of $\Cc_{\Pi}'$ for any $\Pi\in\Pp_{\overline{\Gamma}}$ in order to calculate $min_{C'\in\Cc_{\Pi}'}|C'|$.

\begin{next}\nx{Definition}\label{virtualite}
 {\rm Let $\Gamma$ be a clique complex, $\overline{\Gamma}$ a binomial extension of $\Gamma$, $\FF=\{F_1,\dots,F_k\}$ and $\Pi\in\Pp_{\overline{\Gamma}}$. We say that a cycle $C$ of $\Gamma$ is a ${\bf \Pi}${\bf -virtual minimal cycle} if and only if either $C$ is a minimal cycle of $\Gamma$, or $C$ is not a minimal cycle in $\Gamma$, but $C$ satisfies the following properties:
\begin{enumerate}
 \item Any chord of $C$ is not an edge of $(\overline{\Gamma}^{(\Pi)})_1$.
 \item If $e,e'\in E(C)$ are distinct, then $e,e'$ are not in a same facet of $\Gamma$. 
\end{enumerate}
}
\end{next}

\begin{next}\nx{Remark}\label{virtual_pi}
{\rm Let $\Gamma$ be a clique complex, $\overline{\Gamma}$ a binomial extension of $\Gamma$ and $\Pi,\Pi'\in\Pp_{\overline{\Gamma}}$. If\break $e\in (E(G)\setminus E(\overline{G}^{(\Pi)}))$, then $e\in (E(G)\setminus E(\overline{G}^{(\Pi')}))$. So the definition of a $\Pi$-virtual minimal cycle of  $G$ does not depend on $\Pi$, so a $\Pi$-virtual minimal cycle can just be called {\bf virtual minimal cycle}. We denote by $\Cc^{v}$ the family of all the virtual minimal cycles of $\Gamma$.
}
\end{next}

\begin{next}\nx{Proposition}\label{cyc-vir-min}
  Let $\Gamma$ be a clique complex and $\overline{\Gamma}$ be a binomial extension of $\Gamma$. Then, for any $C\in \Cc^{v}$ there exists $V\subset V(C)$ such that $\Gamma_V$ is a minimal cycle of lenght $\geq 4$.
\end{next}
\proof
Asume that $E(C)=\{\{x_1,x_2\},\{x_2,x_3\},\dots,\{x_s,x_1\}\}$. Let consider the sets $V\subset V(C)$  such that $\Gamma_V$ contains a cycle. If for any of them all the cycles contained in $\Gamma_V$ are not minimal cycles. So, there exist $j\in\{1,\dots,s\}$ such that $\{z_{j},z_{j+2}\}$ is a chord of $C$, i.e. $D=\{z_{j},z_{j+1},z_{j+2}\}$ is a clique of $\Gamma$, so there exists a facet $F$ in $\Gamma$ who contains $D$ in contradiction with the item $2$ of the definition \ref{virtualite}. Thus, there exists $V\subset V(C)$ such that $\Gamma_V$ is a minimal cycle of $\Gamma$.
\qed

\begin{next}\nx{Definition}\label{substitution}
{\rm Let $\Gamma$ be the clique complex and $\overline{\Gamma}$ be a binomial extension of $\Gamma$. Let $\Pi\in\Pp_{\overline{\Gamma}}$, $C$ be a minimal virtual cycle of $\Gamma$ and $e\in E(C)$ and $e=\{x_e,x_e'\}$ be an edge in $E(C)$. We say $e$ is {\bf virtual} if $e\notin E(\overline{\Gamma}^{(\Pi)})$. By the remark \ref{virtual_pi} it does not depend on $\Pi$, i.e. if $e\notin E(\overline{\Gamma}^{(\Pi)})$, then $e\notin E(\overline{\Gamma}^{(\Pi')})$, for any $\Pi'\in\Pp_{\overline{\Gamma}}$. Let $e$ be a virtual edge of $C$, so there exists a unique facet $F_e$ of $\Gamma$ containing $e$. A path $P_e$ from $x_e$ to $x_e'$ in $\overline{G}^{(\Pi)}$ whose vertices are in $\overline{F}_e$ is called a {\bf $\Pi$-local substitution} of $e$ in $\overline{G}^{(\Pi)}$ if there is not an edge in $\overline{G}^{(\Pi)}$ such that it is a chord of $P_e$ and for any $x\in (V(C)\setminus e)$ there is not any edge in $E(\overline{G}^{(\Pi)})$ from $x$ to any vertex in $V(P_e)\setminus e$.}
\end{next}

\begin{next}\nx{Definition}\label{salvar-pellejo}
  {\rm Let $\Gamma$ be a clique complex, $\overline{\Gamma}$ be a binomial extension of $\Gamma$ and $\Cc$ be the family of all virtual minimal cycles of $G$. Let $\Pi\in\Pp_{\overline{\Gamma}}$ and $C$ be a virtual minimal cycle. We set $\Cc_{\Pi}^{\rm ls}(C)$ the family of all minimal cycles $C'$ in $\overline{G}^{(\Pi)}$ obtained from $C$ by replacing all the virtual edges of $C$ by local substitution in $\overline{G}^{(\Pi)}$. We set $\Cc_{\Pi}^{{\rm ls}}=\cup_{C\in\Cc^{v}}\Cc_{\Pi}^{\rm ls}(C)$.}
\end{next}

\begin{next}\nx{Proposition}\atrib{{\bf Going down for minimal cycles}}\label{cicles-ext-per-pi}
Let $\Gamma$ be a clique complex and $\overline{\Gamma}$ be a binomial extension of $\Gamma$. Then, for any minimal cycle $C'$ of lenght $\geq 4$ of $\overline{\Gamma}^{(\Pi)}$, there exists a virtual minimal cycle $C$ of $\Gamma$ and a subset $V'\subset V(\overline{\Gamma})$, such that $C'=(\overline{\Gamma}^{(\Pi)})_{V(C)\cup V'}$.
\end{next}
\proof
Let $C'$ be a minimal cycle of lenght $\geq 4$ of $(\overline{\Gamma}^{(\Pi)})_1$. Set $q:=|C'|$, $$E(C')=\{\{z_1,z_2\},\{z_2,z_3\},\dots,\{z_{q},z_1\}\},$$ and for all $i=1,\dots,q-1$, $e_i=\{z_i,z_{i+1}\}$ and $e_q=\{z_q,z_1\}$. For all $i\in\{1,\dots,q\}$, there exists a facet $F$ of $\Gamma$ such that $e_i\subset \overline{F}$. We assert that there exist at least two facets $F\neq F'$ of $\Gamma$ and two edges $e_i\neq e_j\in E(C')$ such that $e_i\subset \overline{F}$ and $e_{j}\subset \overline{F'}$, otherwise if there exist $F$ facet of $\Gamma$ such that for all edge $e\in E(C)$, $e\in \overline{F}$, so $C$ is a subgraph of $(\overline{\Gamma}^{(\Pi)}_{\overline{F}})_1$ and by the remark \ref{restiction}, $(\overline{\Gamma}^{(\Pi)}_{\overline{F}})_1$ is a chordal graph, thus $C'$ has a chord in $\overline{\Gamma}^{(\Pi)}$ and it is a contradition to the fact that $C'$ is a minimal cycle of $\overline{\Gamma}^{(\Pi)}$. \\

Thus, we have a sequence of facets  $F_{i_1},\dots,F_{i_q}$ in $\Gamma$ such that $e_j\subset \overline{F}_{i_j}$. If $e_{j+1}\subset \overline{F}_{i_j}$, we take $F_{i_{j+1}}=F_{i_j}$. Consequently we can rewrite this sequence like $F_{i_1},\dots,F_{i_{s'}}$, where $q\geq s'\geq 2$ and $F_{i_j}\neq F_{i_{j+1}}$, such that $C'$ is divided by consecutive paths $C_j\subset \overline{F}_{i_j}$ from $z_{i_{j}}$ to $z_{i_{j+1}}$, where $z_1=z_{i_{1}}=z_{i_{s'+1}}$. Moreover, we can suppose that $F_{i_1}\neq F_{i_{s'}}$; since if $F_{i_1}=F_{i_{s'}}$ we can take as initial point of $C'$ the point $z_{i_{s'}}$.
\begin{itemize}
 \item As for all $j\in\{2,\dots,s'+1\}$ $z_{i_j}\in \overline{F}_{j-1}\cap\overline{F}_{j}$, by the definition of $\overline{\Gamma}$, we have that for all $j\in\{1,\dots,s'\}$, $z_{i_j}\in\Gamma$.
 \item We have that $F_{i_{j}}=F_{i_r}$ for any $j,r\in\N$  such that $r\geq j+2$. Otherwise, since $C'$ is a minimal cycle $$\{z_{i_j},z_{i_r}\},\{z_{i_j},z_{i_{r+1}}\},\{z_{i_{j+1}},z_{i_r}\},\{z_{i_{j+1}},z_{i_{r+1}}\}\notin E(\overline{\Gamma}^{(\Pi)})_1 \mbox{ }{\rm(see}\mbox{ }{\rm figure}\mbox{ }{\rm \ref{virtuelcicle})},$$ but these edges are contained in $F_{i_j}$ . 

\begin{figure}[htb]  
  \centering
  \scalebox{0.8}{
  \includegraphics[viewport=220 550 260 650]{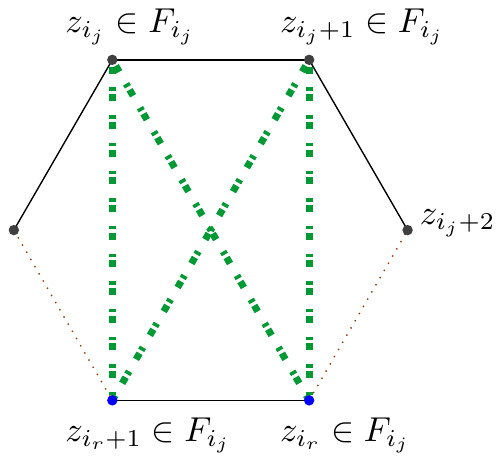}}
  \vspace{3mm}  
  \caption{} 
  \label{virtuelcicle}
\end{figure}

As a result of the definition of $\overline{\Gamma}^{(\Pi)}$, if $e\in (E(\Gamma)\setminus E(\overline{\Gamma}^{(\Pi)}))$ then a vertex of $e$ is in the first line of a matrix $M_{F_{i_j}}$; thus, there exist two elements of $\{z_{i_j},z_{i_{j+1}},z_{i_r},z_{i_{r+1}}\}\subset V(\Gamma)$ which appear in the first line of $M_{F_{i_j}}$. It is a contradiction with the definition of $M_{F_{i_j}}$.
\end{itemize}

Then, we have seen that there exist $s'\geq 3$  and a sequence $F_{i_1},\dots,F_{i_{s'}}$ such that for all\break$j\neq j'\in\{1,\dots,s\}$, $F_{i_j}\neq F_{i_{j'}}$  and $C'$ is divided in consecutive paths $C_j\subset \overline{F}_{i_j}$ from $z_{i_{j}}$ to $z_{i_{j+1}}$, where $z_1=z_{i_{1}}=z_{i_{s'+1}}$. Set $s$ the smallest integer who satisfies these properties.\\
Let $C$ be the cycle defined by $$V(C)=\{z_{i_1},z_{i_2},\dots,z_{i_s}\} \mbox{ }{\rm and}\mbox{ } E(C)=\{\{z_{i_1},z_{i_2}\},\dots,\{z_{i_{s-1}},z_{i_{s}}\},\{z_{i_s},z_{1}\}\}.$$
If $C$ is a minimal cycle, then $C'$ comes from $C$ which is a $\Pi$-virtual minimal cycle by definition \ref{virtualite}. If $C$ is not a minimal cycle, we remark that $C$ is a $\Pi$-virtual minimal cycle:
\begin{enumerate}
 \item Any chord of $C$ is not an edge in $E(\overline{\Gamma}^{(\Pi)})$, since $C'$ is a minimal cycle of $(\overline{\Gamma}^{(\Pi)})$.
 \item By construction of $C$ and by the minimality  of $s$, for all pair $e,e'\in E(C)$, $e$ and $e'$ does not belong to the same facet of $\Gamma$. 
\igno{\begin{itemize}
 \item If $\{z_{i_{j}},z_{i_{j+1}}\},\{,z_{i_{j+1}},z_{i_{j+2}}\}\in E((\overline{\Gamma}^{(\Pi)})_1),$ then $F_{i_1},\dots,F_{i_{j-1}},F,F_{i_{j+2}},\dots,F_{i_s}$ is a sequence of facets of $\Gamma$ with the same characteristiques that $F_{i_1},\dots,F_{i_{s}}$; it is a contradition to the minimality of $s$.
 \item If $\{z_{i_{j}},z_{i_{j+1}}\}\in E((\overline{\Gamma}^{(\Pi)})_1)$ or $\{z_{i_{j+1}},z_{i_{j+2}}\}\in E((\overline{\Gamma}^{(\Pi)})_1)$, but not both. We can suppose that $\{z_{i_{j}},z_{i_{j+1}}\}\in E((\overline{\Gamma}^{(\Pi)})_1)$ and $\{z_{i_{j+1}},z_{i_{j+2}}\}\notin E((\overline{\Gamma}^{(\Pi)})_1)$. Then $\{z_{i_{j+1}},z_{i_{j+2}}\}$ is a proper edge of $\Gamma$, thus $F_{i_{j+1}}=F$. Consequently $F_{i_1},\dots,F_{i_{j-1}},F,$ $F_{i_{j+2}},\dots,F_{i_s}$ is a sequence of facets of $\Gamma$ with the same characteristiques that $F_{i_1},\dots,$ $F_{i_{s}}$; it is a contradiction to the minimality of $s$. 
 \item  $\{z_{i_{j}},z_{i_{j+1}}\},\{,z_{i_{j+1}},z_{i_{j+2}}\}\notin E((\overline{\Gamma}^{(\Pi)})_1).$ Then $\{z_{i_{j}},z_{i_{j+1}}\},\{z_{i_{j+1}},z_{i_{j+2}}\}$ are proper edges of $\Gamma$. Consequently $F_{i_j}=F=F_{i_{j+1}}$, it is not posible because of $F_{i_j}\neq F_{i_{j+1}}$.
\end{itemize}}
\end{enumerate}
Finally, by the last three points $C$ is a $\Pi$-virtual minimal cycle.
\qed\\

\begin{next}\nx{Theorem}\label{sans-linealite}
 Let $\Gamma$ be a clique complex and $\overline{\Gamma}$ a binomial extension of $\Gamma$. Then $\Gamma_1$ is a chordal graph if and only if $\B_{\overline{\Gamma}}$ is $2$-linear.
\end{next}
\proof
\begin{itemize}
 \item [$\Rightarrow)$] As $\Gamma_1$ is a chordal graph, $\Gamma$ does not have minimal cycles. So $(\overline{\Gamma}^{(\Pi)})$ does not too, otherwise, by the proposition \ref{cicles-ext-per-pi} $\Gamma$ would have a $\Pi$-virtual minimal cycle and by the proposition \ref{cyc-vir-min} there exists $V\subset V(C)$ such that $\Gamma_V$ is a minimal cycle of $\Gamma$, but it is not possible, thus $(\overline{\Gamma}^{(\Pi)})$ does not have minimal cycles and by the proposition \ref{betti-p} we obtain that $I_{\overline{\Gamma}^{(\Pi)}}$ is $2$-linear. Then, applying the proposition \ref{igualdad-ini}, we conclude that $\B_{\overline{\Gamma}}$ is $2$-linear.
\item [$\Leftarrow)$] Assume that $\Gamma_1$ is not chordal. So there exist a minimal cycle $C$ of $\Gamma$ of lenght $>3$. As $\overline{\Gamma}$ is an extension of $\Gamma$, all minimal cycle of $G$ is a minimal cycle of $\overline{\Gamma}$, in particular $C$ is a minimal cycle of  $\overline{\Gamma}$, then by the proposition \ref{betti-p},  $\I_{\overline{\Gamma}}$ is not $2$-linear. Moreover by the proposition \ref{primary-decomposition} $$\B_{\overline{\Gamma}}=\bigcap_{{F}\in\FF}{\cal J}_F,$$ where  ${\cal J}_F=(\I_F,(V(\overline{\Gamma})\setminus\overline{F}))$. But the primary decomposition of $\I_{\overline{\Gamma}}$ is: $$\I_{\overline{\Gamma}}=\bigcap_{{F}\in\FF}(V(\overline{\Gamma})\setminus\overline{F}),$$
thus by proposition \ref{2-reg&joint-seq} there does not exist a permutation to the sequence of the ideals $(V(\overline{\Gamma})\setminus\overline{F})$ such that it is linearly joined, in this way by the remark \ref{joint-sequence-heritage} there does not exist a permutation of the sequence $$\{(\I_F,(V(\overline{\Gamma})\setminus\overline{F}))\}_{{F}\in\FF}$$ such that it is linearly joined and so $\B_{\overline{\Gamma}}$ is not  $2$-linear by theorem \ref{2-linear-prime}. \qed
 \end{itemize}

\begin{next}\nx{Remark}\label{subst-virtual}
 {\rm Thanks to proposition \ref{cicles-ext-per-pi} for any minimal cycle $C'$ of $\overline{G}^{(\Pi)}$ there exist a virtual minimal cycle $C$ of $G$, such that $C'\in \Cc_{\Pi}^{\rm ls}(C)$. So, setting $\Cc'_{\Pi}$ the family of all minimal cycles of $\overline{\Gamma}^{(\Pi)}$, then
 \begin{equation}\label{min-con-sust}
  \Cc'_{\Pi}\subset \Cc_{\Pi}^{{\rm ls}}
 \end{equation}
Nevertheless, it is not true in general that $\Cc_{\Pi}^{{\rm ls}}=\Cc'_{\Pi}$, as we can see in the next example:}
\end{next}

\begin{next}\nx{Example}\label{tristement}
 {\rm Let $G$ be the graph of the figure \ref{remplace-bad} and $\Gamma:=\Gamma(G)$ the clique complex generated by $G$.
 
 \begin{figure}[htb]  
  \centering
  \scalebox{0.6}{
  \includegraphics[viewport=340 460 380 650]{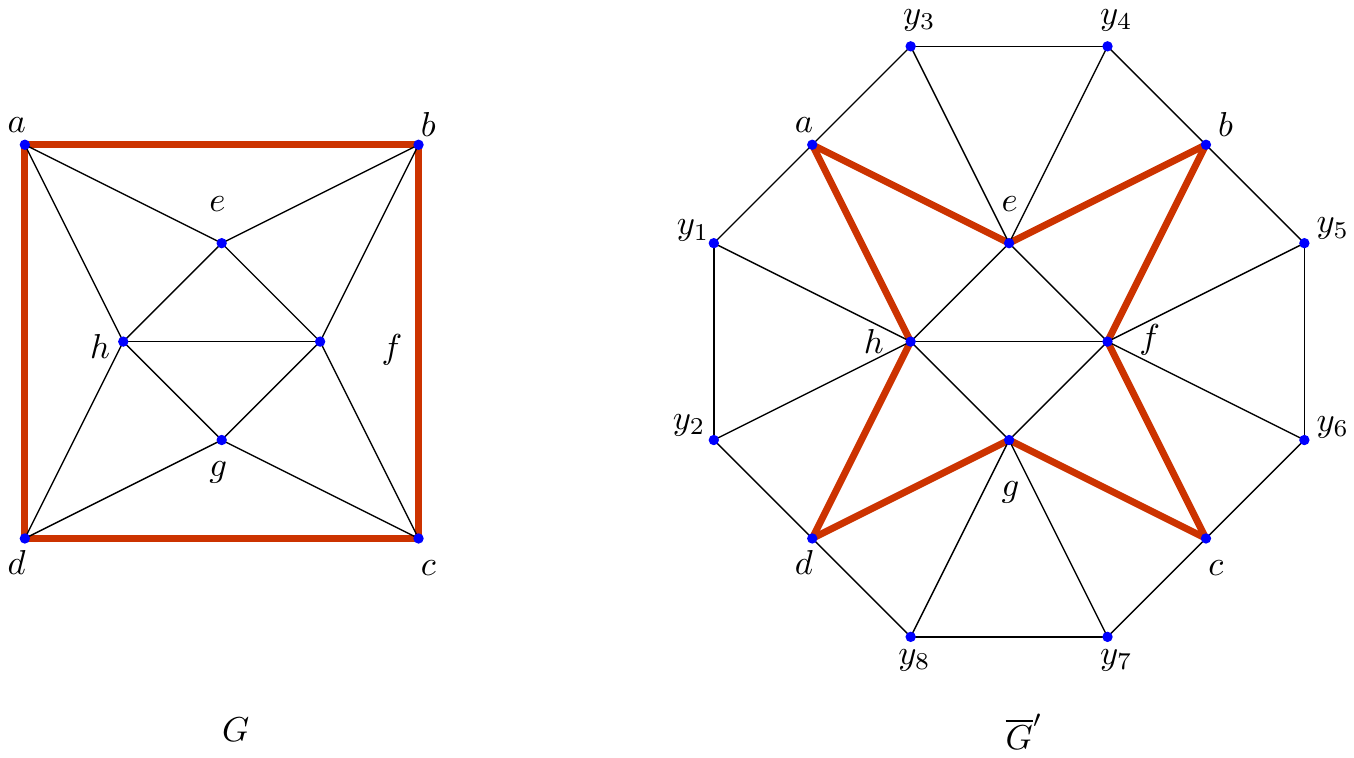}}
  \vspace{3mm}  
  \caption{} 
  \label{remplace-bad}
\end{figure}

 For the facets $F_1=\{a,d,h\}$, $F_2=\{a,b,e\}$, $F_3=\{b,c,f\}$, $F_4=\{c,d,g\}$ we associate the matrices:
 
 \begin{center}
$M_{F_1}=\left(  \begin{array}{lll}
    	a&y_1&y_2\\
	y_1&y_2&d
      \end{array}\right)$,
$M_{F_2}=\left(  \begin{array}{lll}
    	a&y_3&y_4\\
	y_3&y_4&b
      \end{array}\right)$,
\end{center}
\begin{center}
$M_{F_3}=\left(  \begin{array}{lll}
    	b&y_5&y_6\\
	y_5&y_6&c
      \end{array}\right)$ and
$M_{F_4}=\left(  \begin{array}{lll}
    	c&y_7&y_8\\
	y_7&y_8&d
      \end{array}\right)$,
 \end{center}
 
 respectively. $\overline{G}'=(\overline{\Gamma}')_1$ as in the figure \ref{remplace-bad}.\igno{We have that $M_{F_1},M_{F_2},M_{F_3},M_{F_4}$ are admissibly ordered} We have that the cycle $C$ defined by $E(C)=\{\{a,b\},\{b,c\},\break\{c,d\},\{d,a\}\}$ is a minimal cycle of $G$, then $C$ is a virtual cycle of $G$. The paths: $P_1$, $P_2$, $P_3$, $P_4$ defined by: $E(P_1)=\{\{a,e\},\{e,b\}\}$, $E(P_2)=\{\{b,f\},\{f,c\}\}$, $E(P_3)=\{\{c,g\},\{g,d\}\}$, $E(P_4)=\{\{d,h\},\{h,a\}\}$ are local substitutions of $\{a,b\}$, $\{b,c\}$, $\{c,d\}$, $\{d,a\}$, respectively, in $\overline{G}'$ but $C'$ defined by the union of these paths is a cycle of $\overline{G}'$ but it is not a minimal cycle of $\overline{G}'$.}

 \qeds
\end{next}

Now, for any $\Pi\in\Pp_{\overline{\Gamma}}$, we would like to calculate ${\rm min}_{C'\in\Cc_{\Pi}^{{\rm ls}}}(|C'|)$.
The next proposition follows immediately from the inequalities \ref{max-eq}, \ref{max-eq2} and the remark \ref{subst-virtual}.
\begin{next}\nx{Proposition}\label{first-borne}
 Let $\Gamma$ be a clique complex, $\overline{\Gamma}$ a binomial extension of $\Gamma$ and $\FF=\{F\mbox{ }{\rm facet}\mbox{ }{\rm of}\mbox{ }\Gamma: F\neq\overline{F}\}$ is endowed with an admissible order. Then $$p_2(\B_{\overline{\Gamma}})\geq {\rm max}_{\Pi\in\Pp}({\rm min}_{C'\in\Cc'_{\Pi}}(|C'|))\geq {\rm max}_{\Pi\in\Pp}({\rm min}_{C'\in\Cc_{\Pi}^{{\rm ls}}}(|C'|)).$$
\qed
\end{next}

\begin{next}\nx{Definition}\label{replacement}
 {\rm Let $G$ be a graph, $\Gamma=\Gamma(G)$ the clique complex generated by $G$ and $\overline{\Gamma}$ a binomial extension of $\Gamma$. Let $\Pi\in\Pp_{\overline{\Gamma}}$ and $e$ a virtual edge of $G$. A local substitution of $e$ is called a replacement of $e$ in $\overline{G}^{(\Pi)}$ if the lenght of the local substitution is minimal. The lenght of a replacement of $e$ is unique and it is denoted by $t_{\Pi}(e)$. If $e$ is not virtual we set $t_{\Pi}(e):=1$ and in this case by notation abuse we will say that the path $P_e=\{e\}$ is a replacement of $e$.}
\end{next}

\begin{next}\nx{Notation}\label{conj-ver-matrix}
	Let $F$ be a facet of $\Gamma$ such that $F\neq\overline{F}$. We denote by $V_{\Gamma}(M_F)$ the set of all the vertices in $F$ that are in $M_F$. %On the other hand, let $C$ be a minimal virtual cycle of $\Gamma$, we set $\FF_{C}$ the set of all facets which contain a virtual edge of $C$.
\end{next}

\begin{next}\nx{Proposition}\label{existance-borne-remplacement}
Let $\Gamma$ be a clique complex and $\overline{\Gamma}$ a binomial extension of $\Gamma$. Let $\Pi\in\Pp_{\Gamma}$. Let $C\in\Cc^{v}$ and $e\in E(\Gamma)$ be a virtual edge in $E(C)$, such that $e=\{x_0^{(F)},x_k^{(F)}\}$ for some facet of $F$ with $Y^{F}\neq \emptyset$. Then: 
\begin{enumerate}
 \item If there exist $x\in (F\setminus V_{\Gamma}(M_F))$ such that for any vertex $x'\in V(C)$, $\{x,x'\}\notin E(\Gamma)$, then $t_{\Pi}(e)=2$.
 \item Else, let $J_{F}^{{\rm ls}}=\{k\} \cup \{j\in\{1,\dots,m_{F}\}:\forall x'\in (V(C)\setminus e),\mbox{ }\{x_j^{(F)},x'\}\notin E(\Gamma)\}$:
 \begin{enumerate}
		\item If ${\rm min}_{j\in J_{F}^{{\rm ls}}}(|Y_j^{(F)}|)<|Y_k^{(F)}|$, then $t_{\Pi}(e)\leq {\rm min}_{1\leq j\leq m_F}(|Y_j^{(F)}|)+2$.
		\item If ${\rm min}_{j\in J_{F}^{{\rm ls}}}(|Y_j^{(F)}|)= |Y_k^{(F)}|$, then
  $t_{\Pi}(e)\leq {\rm min}_{1\leq j\leq m_F}(|Y_j^{(F)}|)+1$.
	\end{enumerate}
\end{enumerate}
\end{next}
\proof 
\begin{enumerate}
 \item If  there exist $x\in F\setminus V_{\Gamma}(M_F)\neq \emptyset$ such that for any vertex $x'\in V(C)$, $\{x,x'\}\notin E(\Gamma)$, then $\{x_0^{F},x\},\{x,x_k^{F}\}\in E(\overline{G}^{(\Pi)})$ and $\{x_0^{F},x\},\{x,x_k^{F}\}$ is a replacement of $e$ in $\overline{G}^{(\Pi)}$ and $t_{\Pi}(e)=2$.
 \item We can suppose that $M_F$ has one of the following representations: either
$$M_{F}=\left(  \begin{array}{llll|lll|l|lll}
    	x_{0}^{(F)}&y_{1 1}^{(F)}&	\dots&	y_{1 n_{1}}^{(F)}&y_{2 1}^{(F)}&\dots&y_{2 n_{{2}}}^{(F)}&\dots&	y_{m_{F} 1}^{(F)}& \dots&	y_{m_{F} n_{m_{F}}}^{(F)}\\
	y_{1 1}^{(F)}&y_{1 2}^{(F)}&	\dots&	x_{1}^{(F)}&y_{2 2}^{(F)}&	\dots&	x_{2}^{(F)}&\dots&	y_{m_F 2}^{(F)}& \dots&	x_{m_{F}}^{(F)}
      \end{array}\right),\mbox{ }{\rm or}$$ 
$$M_{F}=\left(  \begin{array}{l|lll|l|lll}
    	x_{0}^{(F)}&y_{2 1}^{(F)}&\dots&y_{2 n_{{2}}}^{(F)}&\dots&	y_{m_{F} 1}^{(F)}& \dots&	y_{m_{F} n_{m_{F}}}^{(F)}\\
	x_{1}^{(F)}&y_{2 2}^{(F)}&	\dots&	x_{2}^{(F)}&\dots&	y_{m_F 2}^{(F)}& \dots&	x_{m_{F}}^{(F)}
      \end{array}\right),$$ 
The following argument is valid for both representations of the matrices $M_{F}$. From now to the end of the proof we write $z$ instead of $z^{(F)}$.\\

As $x_0,y_{j 1}$, with $j\in\{2,\dots,m_F\}$ are located only in the first line of $M_F$ and $y_{1,1}$ appears only in the same line and the same column that $x_0$ in $M_F$, we have that for all $j\in\{1,\dots,m_F\}$, $\{x_0,y_{j 1}\}$ is always an edge of $\overline{G}^{(\Pi)}$. As  $\Pi\in\Pp$ (see definition \ref{permutation-permise}) for all $j\in\{1,\dots,m_F\}$ and $i\in\{1,\dots,n_j\}$, $y_{j i}y_{j i+1}$ neither $y_{j n_j}x_j$ are not  diagonals of $\Pi(M_F)$ from the top to the bottom from left to right, so $\{y_{j i},y_{j i+1}\},\{y_{j n_j},x_j\}\in E(\overline{G}^{(\Pi)})$. Finally, since $x_j$ for all $j\in\{1,\dots, m_F\}$ is only in the second line of $M_F$, we have that for all $j\in\{1,\dots,k-1,k+1,\dots,m_F\}$, $\{x_k,x_j\}\in \overline{G}^{(\Pi)}$ and $\{x_0,x_1\}\in \overline{G}^{(\Pi)}$, if $Y_1^{(F)}=\emptyset$; but if $j\in(\{1,\dots,k-1,k+1,\dots,m_F\}\setminus J_{F}^{{\rm ls}})$, $x_j$ is not a vertex of a local substitution of $e$. Then for all $j\in J_{F}^{{\rm ls}}$ we can build  a local substitution $P_j$ from $x_e$ to $x_e'$:
	\begin{itemize}
	 \item $E(P_j)=\{\{x_0,y_{i1}\},\{y_{i1},y_{i2}\},\dots,\{y_{in_i},x_j\},\{x_j,x_k\}\}$ for all $j\neq k$, and
	 \item $E(P_k)=\{\{x_0,y_{k1}\},\{y_{k1},y_{k2}\},\dots,\{y_{k n_k},x_k\}\}$, if $j=k$.
	\end{itemize}
Thus:
	\begin{enumerate}
		\item  If ${\rm min}_{j\in J_{F}^{{\rm ls}}}(|Y_j^{(F)}|)< |Y_k^{(F)}|$, then there exist $j\in J_{F}^{{\rm ls}}$, such that $P_j$ is a local substitution of $e$ and $|P_j|\leq |P_{j'}|$, for all $j'\in J_{F}^{{\rm ls}}$  so $t_{\Pi}(e)\leq |P_j|= {\rm min}_{j\in J_{F}^{{\rm ls}}}(|Y_j^{(F)}|)+2$.
		\item If ${\rm min}_{j\in J_{F}^{{\rm ls}}}(|Y_j^{(F)}|)=| Y_k^{(F)}|$, then $|P_k|\leq |P_j|$ with $j\in J_{F}^{{\rm ls}}$ and $$t_{\Pi}(e)\leq |P_k|=|Y_k^{(F)}|+1={\rm min}_{j\in J_{F}^{{\rm ls}}}(|Y_j^{(F)}|)+1.$$
	\qed
	\end{enumerate}
\end{enumerate}
%\mbox{ }\\
\igno{In this way, setting $\Cc=\{C\mbox{ }{\rm cycles}\mbox{ }{\rm minimal}\mbox{ }\Pi-{\rm virtuel}\mbox{ }{\rm de}\mbox{ }G\}$, on aura que pour toute permutation $\Pi\in\Pp_{\overline{\Gamma}}$:
$$p_2(I_{\overline{\Gamma}^{(\Pi)}})\leq {\rm max}_{\Pi\in\Pp_{\overline{\Gamma}}}\left(p_2(I(\overline{\Gamma}^{(\Pi)}))\right)=p_2(I_{\overline{\Gamma}^{(\Pi^{'})}}).$$}

We recall that $G$ is a graph, $\Gamma=\Gamma(G)$ is the clique complex generated by $G$ and $\overline{\Gamma}$ a binomial extension of $\Gamma$. For each facet $F$ de $\Gamma$ such that $F\neq \overline{F}$, there is a matrix
\begin{center}
$M_F=\left(  \begin{array}{llll|lll|l|lll}
    	x_{0}^{(F)}&y_{{1} 1}^{(F)}&	\dots&	y_{{1} n_{1}}^{(F)}&y_{{2} 1}^{(F)}&\dots&y_{{2} n_{{2}}}^{(F)}&\dots&	y_{{m_F} 1}^{(F)}& \dots&	y_{i_{m_F} n_{{m_F}}}^{(F)}\\
	y_{{1} 1}^{(F)}&y_{1 2}^{(F)}&	\dots&	x_{1}^{(F)}&y_{2 2}^{(F)}&	\dots&	x_{{2}}^{(F)}&\dots&	y_{{m_F} 2}^{(F)}& \dots&	x_{{m_F}}^{(F)}
      \end{array}\right).$
\end{center}
Actually we can write $M_F$ like:
$$M_F=(B_0^{(F)m  }|B_1^{(F)}|\cdots|B_{m_{F}}^{(F)}),$$
where \begin{center}
$B_0^{(F)}=\left(\begin{array}{l}
 x_{i_0}^{(F)}\\
y_{i_{1} 1}^{(F)}\\
\end{array}\right)$ and
$B_j^{(F)}=\left(	\begin{array}{lll}
		y_{i_{j} 1}^{(F)}&	\dots&	y_{i_{j} n_{j_1}}^{(F)}\\
	y_{i_{j} 2}^{(F)}&	\dots&	x_{i_{j}}^{(F)}
		\end{array}\right),$
\end{center}
for all $j=1,\dots,{m_{F}}$.

\begin{next}\nx{Definition}\label{permutation-chidisima}
 {\rm Let $\Gamma$ be a clique complex, $\overline{\Gamma}$ a binomial extension of $\Gamma$ and $F$ a facet of $\Gamma$ different to $\overline{F}$. Let $s={\rm max}\{n_1,\dots,n_{m_F}\}$, we denote by $\Pi^{*}_{F}(M_F)$ the column permutation of the matrix $M_F$ such that:
	$$\Pi^{*}_{F}(M_F)=(B_0^{F}|B_1'^{(F)}|\cdots|B_{s}'^{(F)}),$$
where $j\in\{1,\dots,m_{F}\}$, $B_j'^{(F)}$ is the matrix which columns are all the $j$-th columns of the blocks $B_1^{(F)},\cdots,B_{m_{F}}^{(F)}$ ordered by the standard order. Let $\FF=\{F \mbox{ }{\rm facet}\mbox{ }{\rm of}\mbox{ } \Gamma:\overline{F}\neq F\}=\{F_1,\dots,F_k\}$. We set $\Pi^{*}=(\Pi^{*}_{F_1},\dots,\Pi^{*}_{F_k})$.}
% \item En soyant $\{x_{i,0},x_{i,j}\}$ une ar\^ete prope de $F$
\end{next}

The next propositions is clear from the definition \ref{permutation-chidisima}.
\begin{next}\nx{Proposition}\label{exemples-per-premier-ordre}
 Let $\Gamma$ be a clique complex, $\overline{\Gamma}$ a binomial extension of $\Gamma$ and $F$ a facet of $\Gamma$ such that $F\neq \overline{F}$. Then $\Pi^{*}_{F}$ is an admissible permutation (see definition \ref{permutation-permise}) and consequently $\Pi^{*}\in\Pp_{\overline{\Gamma}}$.
\qed
\end{next}

\begin{next}\nx{Proposition}\label{borne-max}
Let $\Gamma$ be a clique complex and $\overline{\Gamma}$ a binomial extension of $\Gamma$. Let $e$ be a virtual edge of $G$ in $\overline{G}^{(\Pi^{*})}$, such that $e=\{x_0^{(F)},x_k^{(F)}\}$ for some facet $F$ whose $Y_k^{F}\neq \emptyset$. Then:
\begin{enumerate}
  \item If there exist $x\in (F\setminus V_{\Gamma}(M_F))$ such that for any vertex $x'\in V(C)$, $\{x,x'\}\notin E(\Gamma)$, then $t_{\Pi}(e)=2$.
 \item Else, let $J_{F}^{{\rm ls}}=\{k\} \cup \{j\in\{1,\dots,m_{F}\}:\forall x'\in (V(C)\setminus e),\mbox{ }\{x_j^{(F)},x'\}\notin E(\Gamma)\}$:
 \begin{enumerate}
		\item If ${\rm min}_{j\in J_{F}^{{\rm ls}}}(|Y_j^{(F)}|)<|Y_k^{(F)}|$, then  $t_{{\Pi^*}}(e)= {\rm min}_{1\leq j\leq m_F}(|Y_j^{(F)}|)+2$.
		\item If ${\rm min}_{j\in J_{F}^{{\rm ls}}}(|Y_j^{(F)}|)= |Y_k^{(F)}|$, then
  $t_{{\Pi^*}}(e)={\rm min}_{1\leq j\leq m_F}(|Y_j^{(F)}|)+1$.
	\end{enumerate}
\end{enumerate}
\end{next}
\proof
From the proposition \ref{existance-borne-remplacement} we have that the first statement is true. We should only prove the second statement
\begin{enumerate}
 \item [$2.$] 
 \begin{enumerate}
  \item By the definition of $\Pi^{*}_{F}(M_F)$ (see \ref{complexe-sans-diagonal}) we can describe all edges in $\overline{\Gamma}^{(\Pi^*)}$ inside $\overline{F}$, i.e.
   \begin{itemize}
    \item $\{x_0,y_{j1}\}$ with $j\in\{1,\dots,m_F\}$; $\{x_0,x_1\}$ if $Y_1^{(F)}=\emptyset$.
    \item For any $j,t\in\{1,\dots, m_F\}$, $i\in\{1,\dots,n_j\}$, $s\in\{1,\dots,n_t\}$ and $i\geq s$, $\{y_{ji},y_{ts}\}\in \overline{\Gamma}^{(\Pi^*)}$ if and only if 
    \begin{enumerate}
      \item Either $i=s$ and $j\neq t$;
      \item or $i=s+1$ and $t\geq j$;
    \end{enumerate}
     \item For any $j,t\in\{1,\dots, m_F\}$, $i\in\{1,\dots,n_j\}$ and $\varphi(t)$ the number in $\N$ such that $x_t$ appears in the block $\varphi(t)$ of $\Pi^{*}(M_F)$; $\{y_{ji},x_t\}\in\overline{\Gamma}^{(\Pi^*)}$ if and only if $i\geq\eta_e$ and
      \begin{enumerate}
	\item either $\varphi(t)<i$;
	\item or $\varphi(t)=i$ and $j\geq t$.
      \end{enumerate}
    \igno{$\left(\begin{array}{l}y_{ji}\\ 
      z
    \end{array}\right)$ is a column of $\Pi^{*}(M_F)$ or there exists a column $\left(\begin{array}{l}*\\ 
      z
    \end{array}\right)$ preceding a column $\left(\begin{array}{l}y_{ji}\\ 
      *'
    \end{array}\right)$ and there does not exist a column $\left(\begin{array}{l}z\\ 
      *''
    \end{array}\right)$ preceding a column $\left(\begin{array}{l}*'''\\ 
      y_{ji}
    \end{array}\right)$.}
    \item $\{x_j,x_{j'}\}\in E(\overline{\Gamma}^{({\Pi*})})$ for all $1\leq j<j'\leq m_F$
    
    \end{itemize}
    \item We set $\eta_{e}={\rm min}_{j\in J_{F}^{{\rm ls}}}(|Y_j^{(F)}|)$. Let $j'\in J_{F}^{{\rm ls}}$ such that $|Y_{j'}^{(F)}|=\eta_e$, if $|Y_{k}^{(F)}|=\eta_e$, we set $j'=k$. We define $P_e$ as 
 \begin{enumerate}
		\item if $j'\neq k$, $E(P_e)=\{\{x_0,y_{j'1}\},\{y_{j'1},y_{j'2}\},\dots,\{y_{j'\eta_e},x_{j'}\},\{x_{j'},x_k^{(F)}\}\}$; note that\break $|P_e|=\eta_e+2$.
		\item  If $j'=k$, $E(P_e)=\{\{x_e,y_{k1}\},\{y_{k1},y_{k2}\},\dots,\{y_{k\eta_e-1},y_{k \eta_e}\},\{y_{k n_e},x_k^{(F)}\}\}$; note that\break $|P_e|=\eta_e+1$.
 \end{enumerate}
 By $(a)$ $P_e$ is a path in $\overline{\Gamma}^{(\Pi^*)}$ and does not have chords in $F$ because $\Pi^{*}$ is an admissible order. Moreover, since $j'\in J_{F}^{{\rm ls}}$, there is not any edge from $V(P_e)\setminus e$ to a point $x\in V(C)\setminus e$. Thus $P_e$ is a local substitution of $e$.

  \item We are going to prove that $P_e$ is a local substitution of $e$ of minimal lenght.\\
  Let $R$ be any local substitution of $e$. It is clear that $R$ must start with an edge $\{x_{0},y_{j_1,1}\}$ for $1\leq j\leq m_F$ or $\{x_0,x_1\}$ if $1\in J_{F}^{{\rm ls}}$. In the case where $R$ start with $\{x_0,x_1\}$, $\eta_e=0$ and by $(a)$ the stament is true. So we assume that $R$ starts with $\{x_{0},y_{j_1,1}\}$ for some $1\leq j\leq m_F$. By $(a)$ the last edge of $R$ is:
   \begin{itemize}
    \item either $\{y_{ji},x_k\}$ with $i\geq \eta_e$ and either $\varphi(k)<i$, or $\varphi(k)=i$ and $j\geq k$.
    \item Or $\{x_{t},x_k\}$ with $t\in J_{F}^{{\rm ls}}$. In this case there is an edge $\{y_{ji},x_{t}\}$ with $i\geq \eta_e$ and either $\varphi(t)<i$ or $\varphi(t)=i$ and $j\geq t$.
    In any case there is $y_{ji}\in V(R)$, with $i\geq\eta_e$. Moreover, any other edge of $R$ has the form $\{y_{ji},y_{ts}\}$ with either $i=s$ and $j\neq t$; or $i=s+1$ and $t\geq j$. This prove that in order to go to $y_{ji}$ with $i\geq\eta_e$ from $y_{j_1,1}$ we have to go forward or backward block by block, but we can never go forward more than one block by one step. Thus:
     \begin{enumerate}
		\item If ${\rm min}_{j\in J_{F}^{{\rm ls}}}(|Y_j^{(F)}|)<|Y_k^{(F)}|$, then  $|R|\geq i+2\geq\eta_e+2=|P_e|$.
		\item If ${\rm min}_{j\in J_{F}^{{\rm ls}}}(|Y_j^{(F)}|)= |Y_k^{(F)}|$, then $|R|\geq i+1\geq\eta_e+1=|P_e|$.
  	\end{enumerate}
   \end{itemize}
 \end{enumerate}  

 So $t_{{\Pi^*}}(e)=P_e$.
\qed
 \end{enumerate}

\igno{\begin{next}\nx{Proposition}\label{profe-necio}
 Let $\Gamma$ be a clique complex and $\overline{\Gamma}$ a binomial extension of $\Gamma$. Then for all $\Pi\in\Pp_{\Gamma}$, $p_2(I_{\overline{\Gamma}^{(\Pi)}})\leq p_2(I_{\overline{\Gamma}^{(\Pi^{*})}})$. In particular, $${\rm max}_{\Pi\in\Pp_{\overline{\Gamma}}}(p_2(\I_{\overline{\Gamma}^{(\Pi)}}))=p_2(\I_{\overline{\Gamma}^{(\Pi^{*})}})={\rm min}_{C\in\Cc}(C')-3,$$
where $C^{(\Pi^{*})}$ is any cycle of $\Cc^{(\Pi^{*})}(C)$.
\end{next}
\proof
For all $\Pi\in\Pp_{\Gamma}$, let $\Cc$ be the family of all the $\Pi$-virtual cycles of $G$. If $e\in G$ is virtual, by propositions \ref{existance-borne-remplacement} and \ref{borne-max}, we have that  $t_{\Pi}(e)\leq t_{\Pi^{*}}(e)$. Thanks to the proposition \ref{betti-p}, we know that for all $\Pi\in\Pp_{\Gamma}$, $p_2(I_{\overline{\Gamma}^{(\Pi)}})+3$ is equal to the minimal lenght of a minimal cycle of $(\overline{\Gamma}^{(\Pi)})_1=\overline{G}^{(\Pi)}$. By the propositions \ref{importantin} and \ref{borne-max}, for all $\Pi\in\Pp_{\overline{\Gamma}}$,
$$\begin{array}{ll}p_2(I_{\overline{\Gamma}^{(\Pi)}})&={\rm min}_{C\in\Cc}(C_(\Pi)^{r})-3\\
&={\rm min}_{C\in\Cc}(\sum_{e\in E(C)}t_{\Pi}(e))-3\\
&\leq {\rm min}_{C\in\Cc}(\sum_{e\in E(C)}t_{\Pi^{*}}(e))-3\\
&={\rm min}_{C\in\Cc}(C_{\Pi^*}^{r})-3=p_2(I_{\overline{\Gamma}^{(\Pi^{*})}}),
\end{array}
$$
where $C^{(\Pi)}$ is any cycle of $\Cc_{\Pi}^{r}(C)$.
So 
$${\rm max}_{\Pi\in\Pp_{\overline{\Gamma}}}(p_2(I_{\overline{\Gamma}^{(\Pi)}}))=p_2(I_{\overline{\Gamma}^{(\Pi^{*})}})={\rm min}_{C\in\Cc}(C^{(\Pi^*)})-3.$$
\qed}

\begin{next}\nx{Theorem}\label{borne-in-max}
 Let $\Gamma$ be a clique complex and $\overline{\Gamma}$ a binomial extension of $\Gamma$ such that $\FF=\{F \mbox{ }{\rm facet}\mbox{ }{\rm of}\mbox{ } \Gamma:\overline{F}\neq F\}$ is endowed by an admissible order. Then
\begin{enumerate}
 \item $p_2(\B_{\overline{\Gamma}})\geq {\rm min}_{C\in\Cc}({\rm min}_{C'\in\Cc_{\Pi^*}^{\rm ls}(C)}|C'|)-3$.
 \item For any virtual edge $e\in(C)$ and $F$ the unique facet containing it, let $J_{F}^{{\rm ls}}=\{j\in\{1,\dots,m_{F}\}:\forall x'\in (V(C)\setminus e),\mbox{ }\{x_j^{(F)},x'\}\notin E(\Gamma)\}\cup\{k\}$, $\eta_{e}={\rm min}_{j\in J_{F}^{{\rm ls}}}(|Y_j^{(F)}|)$,\\
 $R^{(\Pi^*)}(C)=\{e\in E(C)|e\notin E(\overline{G}^{(\Pi^*)})\}$;\hfill\break
$R_1^{({\Pi^*})}(C)=\{e\in R^{(\Pi^*)}(C)|\exists x\in F\setminus V_{\Gamma}(M_{F}),\mbox{ }\forall x'\in V(C),\mbox{ }\{x,x\}\notin E(\overline{\Gamma})\};$\hfill\break$R_2^{(\Pi^*)}(C)=\{e\in R^{(\Pi^*)}(C)\setminus R_1^{(\Pi^*)}(C)|\exists k\in\{1,\dots,m_{F}\},\mbox{ }e\cap B_k^{(F)}\neq \emptyset,\mbox{ }{\rm and}\mbox{ }|Y_k^{(F)}|>\eta_e\}$;\hfill\break$R_3^{(\Pi^*)}(C)=\{e\in R^{(\Pi^*)}(C)\setminus R_1^{(\Pi^*)}(C)|\exists k\in\{1,\dots,m_{F}\},\mbox{ }e\cap B_k^{(F)}\neq \emptyset,\mbox{ }{\rm and}\mbox{ }|Y_k^{(F)}|=\eta_e\}$.
Then, for any $C\in\Cc$ and any $C^{(\Pi^*)}\in\Cc^{(\Pi^{x})}(C)$, $$C^{(\Pi^*)}=|C|+|R_1^{({\Pi^*})}(C)|+|R_2^{({\Pi^*})}(C)|+\sum_{
 e\in R_{2}^{({\Pi^*})}(C)\cup R_{3}^{({\Pi^*})}(C) }\eta_e.$$
\end{enumerate}
\end{next}
\proof
\begin{enumerate}
 \item Thanks to proposition \ref{first-borne}
$$p_2(\B_{\overline{\Gamma}})\geq {\rm max}_{\Pi\in\Pp_{\overline{\Gamma}}}({\rm min}_{C'\in\Cc_{\Pi}^{\rm ls}}|C'|)-3,$$
%p_2(I_{\overline{\Gamma}^{(\Pi)}}))=p_2(I(\overline{\Gamma}^{(\Pi^{*})})= {\rm min}_{C\in\Cc}(\sum_{e\in C} t_{(\Pi^*)}(e))-3={\rm min}_{C\in\Cc}(C^{(\Pi^*)})-3.
Let $\Pi\in\Pp$. By definition of $\Cc_{\Pi}^{\rm ls}$ (c.f., \ref{salvar-pellejo}) and proposition \ref{borne-max}

   $${\rm min}_{C'\in\Cc_{\Pi^*}^{\rm ls}}|C'|={\rm min}_{C\in\Cc}({\rm min}_{C'\in\Cc_{\Pi}^{{\rm ls}}(C)}|C'|).$$
And\\
$\begin{array}{ll}
{\rm min}_{C'\in\Cc_{\Pi}^{{\rm ls}}(C)}|C'|&={\rm min}_{C\in\Cc}(\sum_{e\in E(C)}t_{\Pi}(e))\\
&\leq {\rm min}_{C\in\Cc}(\sum_{e\in E(C)}t_{\Pi^{*}}(e))\\
&={\rm min}_{C'\in\Cc_{\Pi^{*}}^{{\rm ls}}(C)}|C'|.
  \end{array}$\\
  
So ${\rm max}_{\Pi\in\Pp_{\overline{\Gamma}}}({\rm min}_{C'\in\Cc_{\Pi}^{\rm ls}}|C'|)={\rm min}_{C'\in\Cc_{\Pi^{*}}^{{\rm ls}}(C)}|C'|$ and we conclude that
$$p_2(\B_{\overline{\Gamma}})\geq {\rm min}_{C\in\Cc}({\rm min}_{C'\in\Cc_{\Pi^*}^{\rm ls}(C)}|C'|)-3.$$
 \item Let $C\in\Cc$ and $C^{(\Pi^*)}\in\Cc^{(\Pi^{x})}(C)$, using the proposition \ref{borne-max} we have that
{\small $$\begin{array}{ll}
	\sum_{e\in C} t_{({\Pi^*})}(e)&=\sum_{e\in E(C)-R^{({\Pi^*})}(C)}t_{{\Pi^*}}(e)+\sum_{e\in R_1^{({\Pi^*})}(C)}t(e)+\sum_{e\in R_2^{({\Pi*})}(C)}t_{{\Pi^*}}(e)+\sum_{e\in R_3^{({\Pi^*})}(C)}t_{{\Pi^*}}(e)\\
	\\
	&=\sum_{e\in E(C)-R^{({\Pi^*})}(C)}1+\sum_{e\in R_1^{({\Pi^*})}(C)}2+\sum_{e\in R_2^{({\Pi^*})}(C)}(\eta_e+2)+\sum_{e\in R_3^{({\Pi^*})}(C)}(\eta_e+1)\\
	\\
	&=\sum_{e\in E(C)}1+\sum_{e\in R_1^{({\Pi^*})}(C)}1+\sum_{e\in R_2^{({\Pi^*})}(C)}1+\sum_{e\in R_2^{({\Pi^*})}(C)\cup R_3^{({\Pi^*})}(C)}(\eta_e)\\ 
	\\
	&=|C|+|R_1(C)|+|R_2(C)|+\sum_{e\in R_2^{({\Pi^*})}(C)\cup R_3^{({\Pi^*})}(C)}(\eta_e).
  \end{array}$$}
\end{enumerate}
\qed
%------------------------
%------------------------
\section{Upper bounds of $p_2(\B_{\overline{\Gamma}})$}

\begin{next}\nx{Proposition}\label{poligon}
	Let $C_n$ be the cycle of $n$-vertices $x_1,\dots,x_n$ and $\Gamma$ the clique complex generated by $C_n$. For each edge $\{x_k,x_{k+1}\}$ where $k=1,\dots,n$ ($x_{n+1}=x_1$) we associate a set $Y_{k}$ of cardinality $s_{i}\in\N$ and we set $s=\sum_{j=1}^m s_{i}$. Then $\beta_{i,j}(\B_{\overline{\Gamma}})=\beta_{i,j}(I_{\overline{\Gamma}'})=(d+s)\frac{i}{d+s-i-1}{d+s-2\choose i+1}$ and $p_2(\B_{\overline{\Gamma}})=p_2(I_{\overline{\Gamma}{'}})=n+s-3$. In particular $\si/\B_{\Gamma}$ is Gorenstein.
\end{next}
\proof
By remark \ref{ciclillo}, $\FF=\{F\mbox{ }{\rm facet}\mbox{ }{\rm of}\mbox{ } C_n: \overline{F}\neq F\}$ is endowed by an admissible ordered. So, by theorem \ref{preliminar-in} $\rm{in}(\B_{\overline{\Gamma}})=I_{\overline{\Gamma}'}$. By definition \ref{complexe-sans-diagonal} we have that $\overline{\Gamma}'$ is a $n+s$ cycle. By Hochster's formula we know that the free minimal of resolution $\si/I_{\overline{\Gamma}{'}}$ has the following shape:
	$$0\rightarrow \si(-n-s)\rightarrow \si^{\beta_{n+s-3}}(-n-s+2)\rightarrow\dots\rightarrow\si^{\beta_{1}}(-2)\rightarrow \si\rightarrow \si/I_{\overline{\Gamma}{'}}\rightarrow 0,$$
	where for $i=1,\dots,n+s-3$, 
	\begin{center}
$\beta_{i}=\beta_{i,i+1}(\si/I_{\overline{\Gamma}^{'}})=\beta_{i-1,i+1}(I_{\overline{\Gamma}^{'}})$ and $1=\beta_{n+s-2,n+s}(\si/I_{\overline{\Gamma}^{'}})=\beta_{n+s-3,n+s}(I_{\overline{\Gamma}^{'}})$.
	\end{center}
	Then $p=p_2(I_{\overline{\Gamma}^{'}})=n+s-3$ and ${\rm dimproj}(\si/I)=n+s-2$, thus $\si/I$ is Cohen-Macaulay. By the proposition \ref{igualdad-ini} and by the theorem \ref{betti-gradue-in} we deduce that
	\begin{center}
	$\beta_{i,j}(S/\B_{\overline{\Gamma}})=\left\{\begin{array}{ll}
							0&{\rm si}\mbox{ }j\neq i+1\mbox{ }{\rm and}\mbox{ } i\neq n+s-2,\\
							0&{\rm si}\mbox{ }i= n+s-2\mbox{ }{\rm and} j\neq n+s,\\
							\beta_{i,j}(\si/I_{\overline{\Gamma}^{'}})&{\rm si}\mbox{ }j= i+1.
	                                             \end{array}\right.$
	\end{center}
	Thus ${\rm projdim}\mbox{ }\si/\B_{\bar{\Gamma}}\leq n+s-2$.
	Moreover,  ${\cal K}(\si/\B_{\bar{\Gamma}};t)={\cal K}(\si/I_{\overline{\Gamma}{'}};t)$. So the leading terms are equal and:
$$1=\beta_{n+s-2,n+s}(\si/I_{\overline{\Gamma}^{'}})=\sum_{j=0}^{n+s-2}(-1)^{n+s-2-j}\beta_{j,n+s}=\beta_{n+s-2,n+s}.$$
So the free minimal of resolution of $\si/\B_{\Gamma}$ is:
$$0\rightarrow \si(-n-s)\rightarrow \si^{\beta_{n+s-3}}(-n-s+2)\dots\si^{\beta_{1}}(-2)\rightarrow \si\rightarrow \si/B_{\overline{\Gamma}}\rightarrow 0.$$
As $\si/I_{\Gamma'}$ is Cohen-Macaulay and $d$-pure with $d=(0,2,3,\dots,-2,n+s-2,n+s)$, by the Herzog-K\"uhl equations \cite{herzog1}
$$\beta_i=\beta_0(-1)^{i+1}\prod_{k\neq i,k\neq 0}\frac{d_k-d_0}{d_k-d_i},$$
{\rm for $1\leq i\leq d+s-2=-2$. We set $n:=d+s$ and we have that:}
\begin{center}
 $\begin{array}{ll}
\beta_{i-1}(\B_{\overline{\Gamma}})&=\beta_{i-1}(\I_{\overline{\Gamma}})=\beta_{i}=(-1)^{i+1}(\prod_{k=1}^{i-1}\frac{k+1}{k+1-(i+1)})(\prod_{k=i+1}^{n-3}\frac{k+1}{k+1-(i+1)})(\frac{n-2}{n-2-(i+1)})\\ \\
&=(-1)^{i+1}(-1)^{i-1}(\prod_{k=1}^{i-1}\frac{k+1}{i-k})(\prod_{k=i+1}^{n-3}\frac{k+1}{k-i})(\frac{n}{n-(i+1)})\\ \\
&=\frac{(n-2)!/i+1}{(i-1)!(n-3-i)!}\frac{n}{n-(i+1)}=n\frac{i}{n-i-1}\frac{(n-2)!}{\frac{(i+1)!}{i}(n-2-(i+1))!}\\ \\
&=n\frac{i}{n-i-1}\frac{(n-2)!}{(i+1)!(n-2-(i+1))!}=n\frac{i}{n-i-1}{n-2\choose i+1}.
\end{array}$
\end{center}
In particular, we deduce that $p_2(\B_{\overline{\Gamma}})=n+s-3=p_2(I_{\overline{\Gamma}'})$ and $\si/B_{\overline{\Gamma}}$ is Gorenstein.
\qed\\

Before proving our second principal result (cf. \ref{superieur-borne}) we are going to introduce an special graduation for $\si$. By proposition 1.11 (cf. \cite{eisen-sturm}) if $B$ is a binomial ideal not including monomials, so the algebra $\si$ can be graded by a semigroup $\Sigma(B)$ in the following way: let $B\subset \si$ a binomial ideal not including monomials, so $\Sigma(B)=\N^n/\cong_B$ where 
\begin{center}
 $\forall \alpha,\beta\in\N^n$, $\alpha \cong_B \beta$, if there exist a binomial ${\bf x}^\alpha-a{\bf x}^\beta\in B$ with $a\in\k$.
\end{center}
We set the morphism between semigroups ${\rm deg}_{\Sigma_B}:\N^n\rightarrow \Sigma(B)$ defined by ${\rm deg}_{\Sigma(B)}(\alpha)=[\alpha]_{\cong B}$. So $\si=\oplus_{\sigma\in\Sigma(B)} S_{\sigma}$, where $S_{\sigma}=\oplus_{{\rm deg}(\alpha)=\sigma}\k\cdot x^\alpha$, is a $\Sigma(B)$-graduation of $\si$. If $B$ is an toric homogeneous ideal of $\si$, we have that $\si$ is positively $\Sigma(B)$-graded.

\begin{next}\nx{Proposition} \cite[Proposition 1.1]{bruns}\label{formule-bruns-herzog}
Let $B$ be a toric ideal and $I$ be a squarefree monomial ideal of $\si$ such that: for any binomial ${\bf x}^\alpha-{\bf x}^\gamma\in B$, we have that ${\bf x}^\alpha\in I$ if and only if $x^\gamma\in I$. Let $J=B+I$; $R=\si/J$. Setting 

\begin{center}
$\Omega_b=\{F\subset\{1,\dots,n\}:\exists \beta,\gamma\in\N^n,\mbox{ }{\rm deg}_{\Sigma(B)}(\beta)=b,\mbox{ } {\bf x}^\beta-{\bf x}^\gamma x^{F}\in B\}$ and
$\Psi_b=\{F\in\Omega_b:\exists \beta,\gamma\in\N^n,\mbox{ }{\rm deg}_{\Sigma(B)}(\beta)= b,\mbox{ }{\bf x}^\beta-{\bf x}^\gamma x^{F}\in B, \mbox{ }{\bf x}^\gamma\in I\}.$
\end{center}
Then for all $b\in\Sigma(B)$
$\beta_{i,b}(R)={\rm dim}_{\k}\widetilde{{\rm H}}_{i-1}(\Omega_b,\Psi_b(B);\k)$.
\end{next} 

Now, we return to our problem. Let $\Gamma$ be a clique complex with $n$ vertices, $\overline{\Gamma}$ a binomial extension of $\Gamma$, $m=|V(\overline{\Gamma})|$. Moreover we will denote by ${\bf x}=(x_1,\dots,x_n)$ the variable set identifying the vertex set of $\Gamma$, {\bf y} the variable set identifying the vertices in $V(\overline{\Gamma})\setminus V(\Gamma)$ and ${\bf z}=({\bf x},{\bf y})$. So $I_\Gamma$ is an ideal of $\si$ and  $\B_{\overline{\Gamma}}$ is an ideal of $\sis$. We recall: let $F$ be a facet of $\Gamma$ we set $V_{\Gamma}(M_F)$ the set of all the vertices of $\Gamma$ staying in $M_F$, if $M_F$ is defined; otherwise $V_{\Gamma}(M_F)=\emptyset$. Besides $V_{\overline{\Gamma}}(M_F)$ are all the elements of $\overline{F}$ that appear in $M_F$.

\begin{next}\nx{Definition}\label{some-scroll-torique}
{\rm Let ${\bf z_1},\dots,{\bf z_k}$ be variable sets and $I_1\subset\k[{\bf z_1}],\dots,I_k\subset\k[{\bf z_k}]$ be toric ideals such that for all pair $i,j\in\{1,\dots,k\}$, $i\neq j$, $|{\bf z_i}\cap{\bf z_j}|\leq 1$. We define the {\bf graph of the sequence of toric ideals} $I_1,\dots,I_k$, denoted by $G(I_1,\dots,I_k)$,as the graph whose vertice set is $I_1,\dots,I_k$ and $\{I_i,I_j\}$ is an edge if $I_i$ and $I_j$ has a common variable.}
\end{next}

\begin{next}\nx{Theorem}\label{toriquin}\atrib{\cite{hernan}}
Let ${\bf z_1},\dots,{\bf z_k}$ be variable sets and $I_1\subset\k[{\bf z_1}],\dots,I_k\subset\k[{\bf z_k}]$ be toric ideals such that for any pair $i,j\in\{1,\dots,k\}$, $i\neq j$, $|{\bf z_i}\cap{\bf z_j}|\leq 1$. If for any connex component $G(I_1,\dots,I_k)$ is a tree, then $J:=\sum_{i=1}^k I_i$ is a homogeneous toric ideal and setting $r$ as the number of connexes components of $G(I_1,\dots,I_k)$ we have that:
$${\rm dim}\mbox{ }(\si/J)=\sum_{i=1}^k{\rm dim}\mbox{ }(\k[{\bf z_i}]/I_i)+r-k+1,$$
where ${\bf x}=\cup_{i=1}^k{\bf z_i}$
\end{next}

\begin{next}\nx{Lemma}\label{monom-mineur}
 Let $\Gamma$ be a clique complex, $\overline{\Gamma}$ a binomial extension of $\Gamma$ such that $\J_{\overline{\Gamma}}$ is a toric ideal. Let $F\in \FF=\{F\mbox{ }{\rm facet}\mbox{ }{\rm of}\mbox{ } \Gamma: \overline{F}\neq F\}$. For all $2\times 2$ minor of $M_F$, ${\bf z}^\alpha-{\bf z}^\beta$, and all $\gamma\in \N^{m}$ we have that ${\bf z}^{\gamma+\alpha}\in \I_{\overline{\Gamma}}$ if and only if ${\bf z}^{\gamma+\beta}\in \I_{\overline{\Gamma}}$. In particular if ${\bf z}^{\gamma+\alpha}\in \I_{\overline{\Gamma}}$, then ${\bf z}^\gamma({\bf z}^\alpha-{\bf z}^\beta)\in \I_{\overline{\Gamma}}$.
\end{next}
\proof To show that ${\bf z}^{\gamma+\alpha}\in \I_{\overline{\Gamma}}$ if and only if ${\bf z}^{\gamma+\beta}\in \I_{\overline{\Gamma}}$, it is enough to show that if ${\bf z}^{\gamma+\alpha}\in \I_{\overline{\Gamma}}$ then ${\bf z}^{\gamma+\beta}\in \I_{\overline{\Gamma}}$, since ${\bf z}^\beta-{\bf z}^\alpha=-({\bf z}^\alpha-{\bf z}^\beta)$ is also a $2\time 2$ minor of $M_F$.\\

As ${\bf z}^\alpha-{\bf z}^\beta$ is a $2\times 2$ minor of $M_F$, $\gamma$ has to be different to $0$, otherwise ${\bf z}^\alpha\in I_{\Gamma}$, but by definition of $M_F$, is not possible. Furthermore, as ${\bf z}^{\gamma+\alpha}\in \I_{\overline{\Gamma}}$, there exist $z_1,z_2\in V(\overline{\Gamma})$ such that $z_1z_2\in \I_{\overline{\Gamma}}$ and $z_1z_2|z^{\gamma+\alpha}$. We have three cases:\\

\begin{enumerate}
 \item Let $z_1|{\bf z}^\alpha$ and $z_2|{\bf z}^\alpha$. Thus $z^\alpha\in I_{\Gamma}$, but we have already said that by definition of $M_F$, it is not possible.\\

 \item Let $z_1|{\bf z}^\alpha$ or $z_2|{\bf z}^\alpha$ but no both in the same time. Then, we can suppose $z_1|{\bf z}^\alpha$ and $z_2\nmid {\bf z}^\alpha$; so $z_2|{\bf z}^\gamma$. As $z_1|{\bf z}^\alpha$, we have that $z_1\in \overline{F}$. Moreover $z_1z_2\in \I_{\overline{\Gamma}}$, thus $z_2\notin \overline{F}$ and for all $F'$ facet of $\Gamma$ such that $z_2\in\overline{F'}$ $z_1\notin\overline{F'}$. As ${\bf z}^\alpha-{\bf z}^\beta$ is a minor $2\times 2$ of $M_F$, there exist $z_1',z_2'\in V_{\overline{\Gamma}}(M_F)\subset \overline{F}$ such that ${\bf z}^\beta=z_1'z_2'$ and we have two cases:
	\begin{enumerate}
		\item $z_1'=z_2'$. Then $z_1'\in Y^{(F)}$ and $z_1'$ is only in $\overline{F}$.\\

		\item $z_1'\neq z_2'$. Then $z_1'\notin \overline{F}'$ or $z_2'\notin \overline{F}'$ for all $F'$ facet de $\Gamma$ such that $z_2\in\overline{F}'$. Otherwise, if there exist a facet $F'$ of $\Gamma$ such that $z_2\in\overline{F}'$ and $z_1',z_2'\in \overline{F}'$, then $\{z_1',z_2'\}\in {\overline{F}\cap \overline{F'}}$, so $\{z_1',z_2'\}\in {F\cap F'}$ and $\{z_1',z_2'\}$ is not a proper edge of $\Gamma$. As ${\bf z}^\alpha-z_1'z_2'$ is a $2\times 2$ minor of $M_F$, we have that $x_0^{F}\in\{z_1',z_2'\}$, thus, by definition of $\B_{\overline{\Gamma}}$, $\{z_1',z_2'\}$ is an edge proper of $\Gamma$. Consequently $z_1'\notin \overline{F'}$ or $z_2'\notin \overline{F'}$ for all $F'$ facet of $\Gamma$ such that $z_2\in\overline{F}'$.\\
	\end{enumerate} 

In this way, for the two cases ther exist $z'\in \overline{F}$ such that $z'|{\bf z}^\beta$ and $z'\notin \overline{F'}$ for all $F'$ facet of $\Gamma$ such that $z_2\in\overline{F'}$. Consequently $z'z_2\in I_\Gamma$. Moreover $z'z_2|{\bf z}^{\gamma+\beta}$, so ${\bf z}^{\gamma+\beta}\in \I_{\overline{\Gamma}}$.\\

 \item Let $z_1\nmid {\bf z}^\alpha$ and $z_2\nmid{\bf z}^\alpha$. Then $z_1z_2|{\bf z}^\gamma$. So $z_1z_2|{\bf z}^{\gamma+\beta}$ an in this way ${\bf z}^{\gamma+\beta}\in \I_{\overline{\Gamma}}$.
\end{enumerate}
\qed

\begin{next}\nx{Proposition}\label{monom-binom}
Let be $G$ a graph, $\Gamma=\Gamma(G)$ its clique complex, $\overline{\Gamma}$ a binomial extension of $\Gamma$ such that $\J_{\overline{\Gamma}}$ is a toric ideal. If ${\bf z}^\alpha-{\bf z}^\beta\in \J_{\overline{\Gamma}}$ such that ${\bf z}^\alpha\in \I_{\overline{\Gamma}}$, then ${\bf z}^\beta\in \I_{\overline{\Gamma}}$.
\end{next}
\proof
As ${\bf z}^\alpha-{\bf z}^\beta\in \J_{\overline{\Gamma}}$, there exist $a_1,\dots,a_n\in\k$, $\gamma_1,\dots,\gamma_n\in\N^m$, $F_1,\dots,F_n$ facets of $\Gamma$ and $b_1,\dots,b_n$ $2\time 2$ minors of $M_{F_1},\dots,M_{F_n}$ respectively, such that ${\bf z}^\alpha-{\bf z}^\beta=\sum_{i=1}^{n}a_i{\bf z}^{\gamma_i}b_i$. Let $${\cal K}:=\{i\in \{1,\dots,n:\mbox{ }{\rm one}\mbox{ }{\rm of}\mbox{ }{\rm the}\mbox{ }{\rm monomials}\mbox{ }{\rm of}\mbox{ }z^{\gamma_i}b_i\in \I_{\overline{\Gamma}}\}.$$ ${\cal K}\neq \emptyset$, since there must exist $i\in\{1,\dots,n\}$ such that $z^\alpha$ is a monomial of ${\bf z}^{\gamma_i}b_i$. By the lemma \ref{monom-mineur}, if one of the monomials of ${\bf z}^{\gamma_i}b_i\in \I_{\overline{\Gamma}}$, then the two monomials of $z^{\gamma_i}b_i$ belong to $\I_{\overline{\Gamma}}$. So $$z^\alpha-z^\beta= \sum_{i\in{\cal K}}a_i{\bf z}^{\gamma_i}b_i+\sum_{i\notin{\cal K}}a_i{\bf z}^{\gamma_i}b_i.$$

\begin{itemize}
 \item[\dtb] If $z^\beta$ appears in a term of $\sum_{i\in{\cal K}}a_i{\bf z}^{\gamma_i}b_i$, then $z^\beta\in I_{\Gamma}$.
 \item[\dtb] If $z^\beta$ does not appear as a term of $\sum_{i\in{\cal K}}a_i{\bf z}^{\gamma_i}b_i$, then $z^\beta$ must to appear as a term of $\sum_{i\notin{\cal K}}a_i{\bf z}^{\gamma_i}b_i$. Since all monomial of the first addition does not appear as a term of the second addition, we have:
 $$z^\alpha-\sum_{i\in{\cal K}}a_i{\bf z}^{\gamma_i}b_i=z^\beta-\sum_{i\notin{\cal K}}a_i{\bf z}^{\gamma_i}b_i=0.$$
 Thus $z^\alpha\in \J_{\overline{\Gamma}}$, but it is a contradiction to the fact that $\J_{\overline{\Gamma}}$ is a toric ideal. So this case is not possible.
 
\qed\\
\end{itemize}

\begin{next}\nx{Lemma}\label{cycle-indispensable}
 Let $G$ be a graph, $\Gamma=\Gamma(G)$ be its clique complex, $\overline{\Gamma}$ be a binomial extension of $\Gamma$ with $m=|V(\overline{\Gamma})|$, such that $\J_{\overline{\Gamma}}$ is a toric ideal. Let $C$ be a virtual minimal cycle  of $G$ satisfying that for any edge $e$ virtual, $e=\{x_0^{(F_e)},x_1^{(F_e)}\}$, where $F_e$ is the only facet of $\Gamma$ containing $e$. Let $\widetilde{C}$ be the cycle obtained from $C$ replacing any virtual edge $e$ by the edges $\{x_0^{(F_e)},y_{11}^{(F_e)}\},\{y_{11}^{(F_e)},y_{12}^{(F_e)}\},\dots,$ $\{y_{1n_1}^{(F_e)},x_{1}^{(F_e)}\}$. Then:
$$\beta_{|\widetilde{C}|-3,{\rm deg}_{\Sigma_{\J_{\overline{\Gamma}}}}(z^{V(\widetilde{C})})}(\B_{\overline{\Gamma}})\geq 1.$$

\end{next}

\proof
We set $b={\rm deg}_{\J_{\overline{\Gamma}}}{\bf z}^{V(\widetilde{C})}$. by the proposition \ref{monom-binom}, if ${\bf z}^\gamma-{\bf z}^\beta\in \J_{\overline{\Gamma}}$ and ${\bf z}^\gamma\in \I_{\overline{\Gamma}}$, then ${\bf z}^\beta\in \I_{\overline{\Gamma}}$. In this way we may use the Bruns-Herzog formula (cf. \ref{formule-bruns-herzog}) to calculate $\beta_{|V(\widetilde{C})|-3,b}$:
$$\beta_{|\widetilde{C}|-3,b}(\B_{\overline{\Gamma}})=\beta_{|\widetilde{C}|-2,b}(\k[{\bf z}]/\B_{\overline{\Gamma}})={\rm dim}_{\k}\mbox{ }\widetilde{{\rm H}}_{|V(\widetilde{C})|-3}(\Omega_b,\Psi_b;\k),$$
where
$$\Omega_b=\{F\subset V(\overline{\Gamma}):\exists \beta,\gamma\in\N^m,\mbox{ } {\rm deg}_{\Sigma(\J_{\overline{\Gamma}})}(\beta)=b,\mbox{ } {\bf z}^\beta-{\bf z}^\gamma {\bf z}^{{\rm supp}(F)}\in \J_{\overline{\Gamma}}\} \mbox{ }{\rm and}$$
$$\Psi_b=\{F\in\Omega_b:\exists \beta,\gamma\in\N^m,\mbox{ }{\rm deg}_{\Sigma(\J_{\overline{\Gamma}})}(\beta)= b,\mbox{ }{\bf z}^\beta-{\bf z}^\gamma {\bf z}^{{\rm supp}(F)}\in \J_{\overline{\Gamma}}, \mbox{ }{\bf z}^\gamma\in \I_{\overline{\Gamma}}\}.$$
\begin{itemize}
 \item[\dtb] Firstly we show that: $\Omega_b=\{F\subset V(\overline{\Gamma}):\exists \gamma\in\N^m,\mbox{ } {\bf z}^\alpha-{\bf z}^\gamma {\bf z}^{{\rm supp}(F)}\in \J_{\overline{\Gamma}}\}$.\\

${\rm deg}_{\Sigma(\J_{\overline{\Gamma}})}(\beta)=b$ if and only if $({\bf z}^\alpha-{\bf z}^{\beta})\in \J_{\overline{\Gamma}}$. Then ${\bf z}^\beta-{\bf z}^\gamma {\bf z}^{{\rm supp}(F)}\in \J_{\overline{\Gamma}}$ if and only if\break $({\bf z}^\alpha-{\bf z}^\gamma {\bf z}^{{\rm supp}(F)})\in \J_{\overline{\Gamma}}$. Consequently  $\Omega_b=\{F\subset V(\overline{\Gamma}):\exists \gamma\in\N^m,\mbox{ } {\bf z}^\alpha-{\bf z}^\gamma {\bf z}^{{\rm supp}(F)}\in \J_{\overline{\Gamma}}\}$.\\

\item[\dtb] We are going to show that for any facet $F$ of $\Gamma$ such that $M_F\neq 0$ we have that any $2\times 2$ minor of $M_F$ $b({\bf z})=z_{i,1}z_{i,2}-z_{i,3}z_{i,4}$ satisfies that: $z_{i,1},z_{i,2}\in V(\widetilde{C})$ if and only if $z_{i,3},z_{i,4}\in V(\widetilde{C})$.

Firstly we recall the definition of $M_F$:
$$M_{F}=\left(  \begin{array}{llll|lll|l|lll}
    	x_{0}^{(F)}&y_{1 1}^{F}&	\dots&	y_{1 n_{1}}^{F}&y_{2 1}^{F}&\dots&y_{2 n_{{2}}}^{F}&\dots&	y_{m_{F} 1}^{F}& \dots&	y_{m_{F} n_{m_{F}}}^{F}\\
	y_{1 1}^{F}&y_{1 2}^{F}&	\dots&	x_{1}^{F}&y_{2 2}^{F}&	\dots&	x_{2}^{F}&\dots&	y_{m_{F} 2}^{F}& \dots&	x_{m_{F}}^{F}
      \end{array}\right);$$ 
and $$\I_{F}=\{[M_{F}]_{1,u}[M_{F}]_{2,v}-[M_{F}]_{1,v}[M_{F}]_{2,u}:1\leq u<v\leq (\sum_{j=1}^{m_{F}}n_j) +1\}.$$
As the restiction of $\widetilde{C}$ in $F$ is $\{x_0^{(F)},y_{11}^{(F)}\},\{y_{11}^{(F)},y_{12}^{(F)}\},\dots,\{y_{1n_1}^{(F)},x_{1}^{(F)}\}$, any minor $b(z)$ of $M_F$ has all its variables in $V(\widetilde{C})$ if and only if $b(z)$ is a minor of the first block of $M_F$; so the second monomial  $b(z)$ have also all its variables in $V(\widetilde{C})$. So we conclude the assertion.\\

\item[\dtb] We are going to show that $\Omega_b=<V(\widetilde{C})>$ (the simplicial complex whose all faces are all the subsets of $V(\widetilde{C})$).\\

Let ${\bf z}^\alpha-{\bf z}^\beta \in \J_{\overline{\Gamma}}$. We suppose that some variable of ${\bf z}^\beta$ is not in $V(\widetilde{C})$. As ${\bf z}^\alpha-{\bf z}^\beta\in \J_{\overline{\Gamma}}$, there exist $a_1,\dots,a_n\in\k$, $\gamma_1,\dots,\gamma_n\in\N^m$, $F_{s_1},\dots,F_{s_n}$ facets of $\Gamma$ and $b_1={\bf z}^{\alpha_1}-{\bf z}^{\beta_1},\dots,\break b_n={\bf z}^{\alpha_n}-{\bf z}^{\beta_n}$ $2\times 2$ minors of the matrices $M_{F_{s_1}},\dots,M_{F_{s_n}}$ respectively, such that\break ${\bf z}^\alpha-{\bf z}^\beta=\sum_{i=1}^{n}a_i{\bf z}^{\gamma_i}b_i$.  Let $${\cal K}:=\{i\in\{1,\dots,n\}:\mbox{ }{\rm a}\mbox{ }{\rm monomial}\mbox{ }{\rm of}\mbox{ }{\bf z}^{\gamma_i}b_i\mbox{ }{\rm has}\mbox{ }{\rm all}\mbox{ }{\rm its}\mbox{ }{\rm variables}\mbox{ }{\rm in}\mbox{ }V(\widetilde{C})\}.$$ 
${\cal K}\neq \emptyset$, since there must exist $i\in\{1,\dots,n\}$ such that ${\bf z}^\alpha$ is a term of any ${\bf z}^{\gamma_i}b_i$. Moreover, we note that if a monomial of ${\bf z}^{\gamma_i}b_i$ has all its variables in $V(\widetilde{C})$, then all the variables of ${\bf z}^{\gamma_i}$ are in $V(\widetilde{C})$ and by the last statement, any variable of the two monomials of $b_i$ are in $V(\widetilde{C})$. Thus, we can write $z^\alpha-z^\beta=\sum_{i\in K}a_i{\bf z}^{\gamma_i}b_i+\sum_{i\notin K}a_i{\bf z}^{\gamma_i}b_i$. And by this remark, if $i\notin K$, any monomial of ${\bf z}^{\gamma_i}b_i$ does not have all its variables in $V(\widetilde{C})$, so:
$$z^\alpha+\sum_{i\in K}a_i{\bf z}^{\gamma_i}b_i=z^\beta+\sum_{i\notin K}a_i{\bf z}^{\gamma_i}b_i=0.$$ Thus $z^\alpha\in \J_{\overline{\Gamma}}$, that is a contradiction to the fact that $\J_{\overline{\Gamma}}$ is a toric ideal. Then all the variables of ${\bf z}^\beta$ are in $V(\widetilde{C})$. And in this we way we can conclude:
$$\begin{array}{ll}
   	\Omega_b&=\{F\subset V(\overline{\Gamma}):\exists \gamma\in\N^m,\mbox{ } {\bf z}^\alpha-{\bf z}^\gamma {\bf z}^F\in \J_{\overline{\Gamma}}\}\\
	&=\{F\subset V(\widetilde{C}):\exists \gamma\in\N^m,\mbox{ } {\bf z}^\alpha-{\bf z}^{\gamma+{\rm supp}(F)}\in \J_{\overline{\Gamma}}\}\\
	&=<V(\widetilde{C})>.
  \end{array}
$$		

\item[\dtb] Let us show that
$\Psi_b=({\overline{\Gamma}}_{V(\widetilde{C})})^{A}$.\\
By the same way that we showed that $\Omega_b=\{F\subset V(\overline{\Gamma}):\exists \gamma\in\N^s,\mbox{ } {\bf z}^\alpha-{\bf z}^\gamma {\bf z}^{{\rm supp}(F)}\in \J_{\overline{\Gamma}}\}$, we can prove that
$$\Psi_b=\{F\in\Omega_b:\exists \gamma\in\N^m, \mbox{ }{\bf z}^\alpha-{\bf z}^\gamma {\bf z}^{{\rm supp}(F)}\in \J_{\overline{\Gamma}} \mbox{ }{\rm and}\mbox{ }{\bf z}^\gamma\in \I_{\overline{\Gamma}}\}.$$
As $\Omega_b=<V(\widetilde{C})>$, we have that 
$$\begin{array}{ll}
   \Psi_b&=\{F\in V(\widetilde{C}):\exists \gamma\in\N^m, \mbox{ }{\bf z}^{{\rm supp}(F)+{\rm supp}(V(\widetilde{C})\setminus F)}-{\bf z}^\gamma {\bf z}^{{\rm supp}(F)}\in \J_{\overline{\Gamma}}\mbox{ }{\rm and}\mbox{ }{\bf z}^\gamma\in \I_{\overline{\Gamma}}\}\\
	&=\{F\in V(\widetilde{C}):\exists \gamma\in\N^m, \mbox{ }{\bf z}^{{\rm supp}(F)}({\bf z}^{{\rm supp}(V(\widetilde{C})\setminus F)}-{\bf z}^\gamma)\in \J_{\overline{\Gamma}} \mbox{ }{\rm and}\mbox{ }{\bf z}^\gamma\in \I_{\overline{\Gamma}}\}.
  \end{array}$$
But $\J_{\overline{\Gamma}}$ is a toric ideal, so $$\J_{\overline{\Gamma}}=(\J_{\overline{\Gamma}}:(z_1\cdots z_m)^\infty)=\{p({\bf z})\in\k[z_1\cdots z_m]: \exists\alpha\in\N^m,\mbox{ } {\bf z}^\alpha p(z)\in \J_{\overline{\Gamma}} \},$$ then
$$\Psi_b=\{F\in V(\widetilde{C}):\exists \gamma\in\N^m, \mbox{ }({\bf z}^{{\rm supp}(V(\widetilde{C})\setminus F)}-{\bf z}^\gamma)\in \J_{\overline{\Gamma}}\mbox{ }{\rm and}\mbox{ }{\bf z}^\gamma\in \I_{\overline{\Gamma}}\}.$$
Besides, by the proposition \ref{monom-binom} we deduce:
$$\begin{array}{ll}
	\Psi_b&=\{F\in V(\widetilde{C}):\exists \gamma\in\N^m, \mbox{ }({\bf z}^{{\rm supp}(V(\widetilde{C})\setminus F)}-{\bf z}^\gamma)\in \J_{\overline{\Gamma}}\mbox{ }{\rm and}\mbox{ } {\bf z}^\gamma,{\bf z}^{{\rm supp}(V(\widetilde{C})\setminus F)} \in \I_{\overline{\Gamma}}\}\\
	&=\{F\in V(\widetilde{C}): {\bf z}^{{\rm supp}(V(\widetilde{C})\setminus F)} \in \I_{\overline{\Gamma}}\}\\
	&=\{F\in V(\widetilde{C}):(V(\widetilde{C})\setminus F) \notin {\overline{\Gamma}}\}\\
	&=\{F\in V(\widetilde{C}):(V(\widetilde{C})\setminus F) \notin {\overline{\Gamma}}_{V(\widetilde{C})}\}\\
	&=({\overline{\Gamma}}_{V(\widetilde{C})})^{A}.
\end{array}$$
\end{itemize}
Using the Alexander duality and the long exact sequence of relative homology we have that

$$\begin{array}{ll}
	\beta_{|\widetilde{C}|-3,b}(\B_{\overline{\Gamma}})&={\rm dim}_{\k}\widetilde{{\rm H}}_{|\widetilde{C}|-3}(\Omega_b,\Psi_b;\k)\\
	&={\rm dim}_{\k}\mbox{ }\widetilde{{\rm H}}_{|\widetilde{C}|-(|\widetilde{C}|-3)-2}((\Psi_b)^{A},(\Omega_b)^{A};\k)\\
	&={\rm dim}_{\k}\mbox{ }\widetilde{{\rm H}}_{1}((\Psi_b)^{A};\k)\\
	&={\rm dim}_{\k}\mbox{ }\widetilde{{\rm H}}_{1}(\overline{\Gamma}_{V(\widetilde{C})};\k).
  \end{array}
$$
\begin{itemize}
 \item[\dtb] We are going to prove that ${\rm dim}_{\k}\mbox{ }\widetilde{{\rm H}}_{1}(\overline{\Gamma}_{V(\widetilde{C})};\k)\geq 1$.\\
	We note that the facets of $\overline{\Gamma}_{V(\widetilde{C})}$ are the subset $\widetilde{F}_e$ for all $e\in E(C)$ defined by: 
$$\widetilde{F}_e=\left\{\begin{array}{ll}
 		e& {\rm if}\mbox{ }e\in E(\widetilde{C}),\\
		e\cup Y_1^{(F_e)}& {\rm if}\mbox{ }e\notin E(\widetilde{C}).
\end{array}\right.
$$
 So $\overline{\Gamma}_{V(\widetilde{C})}=<\widetilde{F}_e:e\in C>$ and is homotopically equivalent to $\Gamma_{V(C)}$, as we can see in the figure \ref{coll-import-dib}, where we set $E(C)=\{\{x_1,x_2\},\{x_2,x_3\},\dots,\{x_n,x_1\}\}$ and
\begin{center}
if $\{x_i,x_{i+1}\}\in E(\widetilde{C})$, $\widetilde{F}_{x_i,x_{i+1}}=\{x_i,x_{i+1},y_{11}^i,\dots,y_{1n_1}^i\}$, where $x_{n+1}=x_1$.
\end{center}

\begin{figure}[h,ht]  
  \centering
  \scalebox{0.80}{
  \includegraphics[viewport=280 550 340 660]{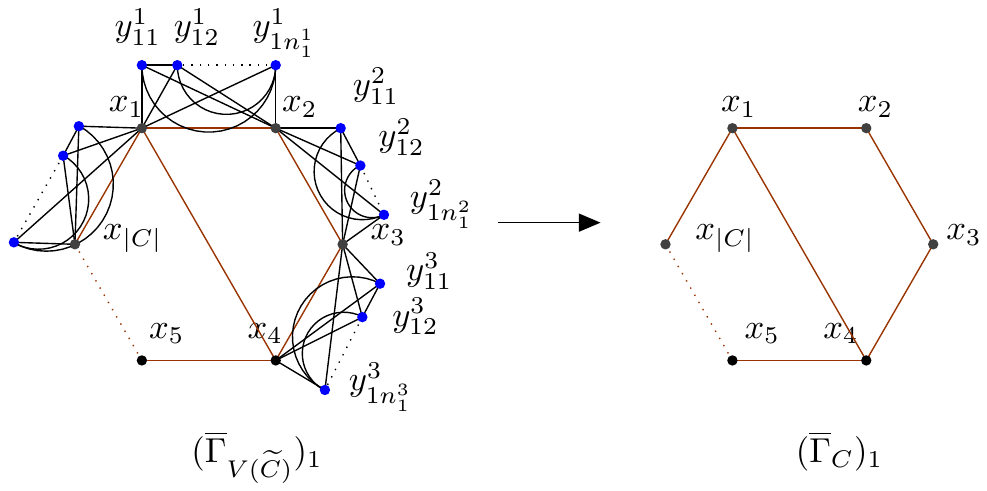}}
  \vspace{3mm}  
  \caption{} 
  \label{coll-import-dib}
\end{figure}

We know that $\widetilde{H}_1(\Gamma_{V(C)};\k)$ is the homology of the complex 
$$\oplus_{\sigma\in (\Gamma_{V(C)})_{(2)}}\k\cdot f_{\sigma}\xrightarrow{\partial_{2}} \oplus_{\sigma\in E(\Gamma_{V(C)})}\k\cdot f_{\sigma}\xrightarrow{\partial_1} \oplus_{\sigma\in {V(C)}}\k\cdot f_{\sigma},$$
where $\Gamma_{V(C)})_{(2)}=\{\sigma\in\Gamma_{V(C)}:|\sigma|=3\}$ and $\partial_i(e_\sigma)=\sum_{i\in\sigma}{\rm sign}(i,\sigma)e_{\sigma\setminus\{i\}},$ for any $i\in\{1,2\}$ and ${\rm sign}(i,\sigma)=(-1)^{r-1}$, if $i$ is the $r$-th element of the set $\sigma\subset\{x_1<\dots<x_n\}$, writen by increassing order.\\

We have that $\partial_1(\sum_{i=1}^{n-1}f_{\{x_i,x_{i+1}\}}-f_{\{x_i,x_{i+1}\}})=0$, so $\sum_{i=1}^{n-1}f_{\{x_i,x_{i+1}\}}-f_{\{x_1,x_{n}\}}\in {\rm ker}(\partial_1)$.\break Now we need to prove that there is not an element $x\in \oplus_{\sigma\in (\Gamma_{V(C)})_2}\k\cdot f_{\sigma}$ such that\break $\partial_2(x)=\sum_{i=1}^{n-1}f_{\{x_i,x_{i+1}\}}-f_{\{x_1,x_{n}\}}$. Otherwise, assume that there exist $x=\sum_{\sigma\in (\Gamma_{V(C)})_{(2)}}a_{\sigma}f_\sigma$ such that $\partial_2(x)= \sum_{i=1}^{n-1}f_{\{x_i,x_{i+1}\}}-f_{\{x_1,x_{n}\}}$. 

We know that $\widetilde{H}_1(\Gamma_{V(C)};\k)$ is the homology of the complex 
$$\oplus_{\sigma\in (\Gamma_{V(C)})_{(2)}}\k\cdot f_{\sigma}\xrightarrow{\partial_{2}} \oplus_{\sigma\in E(\Gamma_{V(C)})}\k\cdot f_{\sigma}\xrightarrow{\partial_1} \oplus_{\sigma\in {V(C)}}\k\cdot f_{\sigma},$$
where $(\Gamma_{V(C)})_{(2)}=\{\sigma\in\Gamma_{V(C)}:|\sigma|=3\}$ and $\partial_i(e_\sigma)=\sum_{i\in\sigma}{\rm sign}(i,\sigma)e_{\sigma\setminus\{i\}},$ for any $i\in\{1,2\}$ and ${\rm sign}(i,\sigma)=(-1)^{r-1}$, if $i$ is the $r$-th element of the set $\sigma\subset\{x_1<\dots<x_n\}$, written by increasing order.\\

We have that $\partial_1(\sum_{i=1}^{n-1}f_{\{x_i,x_{i+1}\}}-f_{\{x_i,x_{i+1}\}})=0$, so $\sum_{i=1}^{n-1}f_{\{x_i,x_{i+1}\}}-f_{\{x_1,x_{n}\}}\in {\rm ker}(\partial_1)$.\break
 Now we need to prove that there is not an element $x\in \oplus_{\sigma\in (\Gamma_{V(C)})_2}\k\cdot f_{\sigma}$ such that
\break $\partial_2(x)=\sum_{i=1}^{n-1}f_{\{x_i,x_{i+1}\}}-f_{\{x_1,x_{n}\}}$. On the contrary, assume that there exist $x=\sum_{\sigma\in (\Gamma_{V(C)})_{(2)}}a_{\sigma}f_\sigma$ such that 
$\partial_2(x)= \sum_{i=1}^{n-1}f_{\{x_i,x_{i+1}\}}-f_{\{x_1,x_{n}\}}$.\\

Let remark that since $C$ is a virtual minimal cycle, any chord of $C$ is not a chord in $\widetilde{C}$; then any  chord $\{A,B\}$ of $C$,  is a proper edge of $C$,
 and if $F$ is the unique facet containing $\{A,B\}$, then either $A=[M_F]_{11}$ or $B=[M_F]_{11}$.\\

Let $e\in E(C)$, since $e$ appears in $\partial_2(x)$  there exist $\sigma=\{x_{i_1},x_{i_2},x_{i_3}\}\in(\Gamma_{V(C)})_{(2)}$ such that $a_\sigma\neq 0$ and  $\{x_{i_1},x_{i_2}\}=e$.
 Neither $\{x_{i_2},x_{i_3}\}\in E(C)$ or $\{x_{i_1},x_{i_3}\}\in E(C)$, by item  2. of definition \ref{virtualite}. 
Thus $\{x_{i_2},x_{i_3}\}$ and $\{x_{i_1},x_{i_3}\}$ are chords of $C$; since  $\{x_{i_2},x_{i_3}\}$ is in the image $\partial_2(\sigma )$ but doesn't appears in $\partial_2(x)$, then
   there exist $\sigma'=\{x_{i_2},x_{i_3},x_{i_4}\} \in (\Gamma_{V(C)})_{(2)}$  with $a_\sigma'\neq 0$. 
By item 1. of definition \ref{virtualite}, $\{x_{i_2},x_{i_3}\}$ is a chord of   $C$, so it is a  proper edge, and  there exist a unique facet $F$ of $\Gamma$ such that
 $\sigma,\sigma'\subset F$ and $|F|\geq 4$. As $\{x_{i_2},x_{i_3}\}$ and $\{x_{i_1},x_{i_3}\}$ are chords of $C$, we should have  $[M_{F}]_{11}=x_{i_3}$. 
 By item 1. of definition \ref{virtualite}, $\{x_{i_4},x_{i_1}\}$ is a chord of $C$ and thus  either $x_{i_4}=[M_{F}]_{11}$ or $x_{i_1}=[M_{F}]_{11}$, but $[M_{F}]_{11}=x_{i_3}$. This is a contradiction. 
Thus, there is not an element $x\in \oplus_{\sigma\in (\Gamma_{V(C)})_2}\k\cdot f_{\sigma}$ such that $\partial_2(x)=\sum_{i=1}^{n-1}f_{\{x_i,x_{i+1}\}}-f_{\{x_1,x_{n}\}}$. So, ${\rm dim}_k\widetilde{H}_1(\Gamma_{V(C)};\k)\geq 1$.
\end{itemize}
Thus $\beta_{|V(\widetilde{C})|-3,b}(\B_{\overline{\Gamma}})={\rm dim}_{\k}\widetilde{{\rm H}}_{|V(\widetilde{C})|-3}(\Omega_b,\Psi_b;\k)={\rm dim}_{\k}\mbox{ }\widetilde{{\rm H}}_{1}(\overline{\Gamma}_{V(\widetilde{C})};\k)\geq 1$.

\qed
\begin{next}\nx{Definition}\label{ay-caray}
We define $\widetilde{\Cc}=\{C \in\Cc^{v}:\forall e\in E(C)\mbox{ }{\rm virtual}\mbox{ } e\in B_1(M_{F_e})\}.$
\end{next}

\begin{next}\nx{Theorem}\label{superieur-borne}
 Let $\Gamma$ be a clique complex and $\overline{\Gamma}$ a binomial extension of $\Gamma$ such that $\J_{\overline{\Gamma}}$ is a toric ideal. If $\widetilde{\Cc}\neq\emptyset$, then $p_2(\B_{\overline{\Gamma}})+3\leq {\rm min}_{C\in\widetilde{\Cc}}(|\widetilde{C}|)={\rm min}_{C\in\widetilde{\Cc}}(|C|+\sum_{e\in R(C)} Y^{F_e}_1)$.
\end{next}
\proof
Let $C\in\widetilde{\Cc}$. Consider $\widetilde{C}$ the cycle defined in lemma \ref{cycle-indispensable}, so by this lemma:
$$\begin{array}{ll}
\beta_{|\widetilde{C}|-3,|\widetilde{C}|}(\B_{\overline{\Gamma}})&=\sum_{h\in\Sigma(\J_{\overline{\Gamma}}),|h|=|\widetilde{C}|}\beta_{|\widetilde{C}|-3,h}(\B_{\overline{\Gamma}})\\
&\geq \beta_{|\widetilde{C}|-3,{\rm deg}_{\Sigma_{\J_{\overline{\Gamma}}}}(z^{{\rm supp}(V(\widetilde{C}))})}(\B_{\overline{\Gamma}})=1,
\end{array}$$
then $p_2(\B_{\overline{\Gamma}})\leq |\widetilde{C}|-3=|C|+\sum_{e\in R(C)} Y^{(F_e)}_1-3$. So $$p_2(\B_{\overline{\Gamma}})+3\leq {\rm min}_{C\in\widetilde{\Cc}}(|\widetilde{C}|)={\rm min}_{C\in\widetilde{\Cc}}(|C|+\sum_{e\in R(C)} Y^{(F_e)}_1).$$

\qed

\begin{next}\nx{Theorem}\label{2pregul2}
 Let $\Gamma$ be a clique complex and $\overline{\Gamma}$ be a binomial extension of $\Gamma$ such that $\FF=\{F \mbox{ }{\rm facet}\mbox{ }{\rm of}\mbox{ } \Gamma:F\neq\overline{F}\}$ is endowed with an admissible order and $\J_{\overline{\Gamma}}$ is a toric ideal. Suppose that $\widetilde{\Cc}=\Cc^{v}$, moreover for any $C\in \widetilde{\Cc}$ and for any edge $e\in E(C)$ with $F_e$ is the only facet that contains $e$, $e$ satisfies one of the next properties:
\begin{enumerate}
 \item either $|Y_1^{(F_e)}|={\rm min}_{1\leq j\leq m_{F_e}}(Y_i^{(F_e)})\geq 2$ and $V_{\Gamma}(M_{F_e})=F_e$, o\`u $J_{F_e}^{{\rm ls}}=\{j\in(\{1,\dots,m_{F_e}\}\setminus\{k\}):\forall x'\in (V(C)\setminus e),\mbox{ }\{x_j^{(F)}\}\cup\{1\}$;
 \item or $|Y_1^{(F_e)}|=1$.
\end{enumerate}
Then, $p_2(\B_{\overline{\Gamma}})={\rm min}_{C\in\Cc^{v}}(|\widetilde{C}|)-3,$ where $\Cc$ is the family of all virtual minimal cycles of $G$.
\end{next}
\proof
%Let $$\widetilde{\Cc}=\{C \mbox{ }{\rm virtual}\mbox{ }{\rm minimal}\mbox{ }{\rm cycle}:\forall e\in E(C)\mbox{ }{\rm virtual}\mbox{ } e\in B_1(M_{F_e})\}.$$
Let $C\in\Cc^{v}$, thanks to definition of $\Cc_{\Pi^{*}}^{\rm ls}(C)$, $\widetilde{C}\in \Cc_{\Pi}^{\rm ls}(C)$. By proposition \ref{borne-max} $|\widetilde{C}|={\rm min}_{C'\in \Cc_{\Pi^{*}}^{\rm ls}(C)}|C'|$. Thus, from the theorems \ref{superieur-borne} and \ref{borne-in-max}
$$p_2(\B_{\overline{\Gamma}})\leq {\rm min}_{C\in\Cc^{v}}(|\widetilde{C}|)-3={\rm min}_{C\in\Cc^{v}}({\rm min}_{C'\in \Cc_{\Pi^{*}}^{\rm ls}(C)}|C'|)-3\geq p_{2}(\B_{\overline{\Gamma}}).$$
So $p_2(\B_{\overline{\Gamma}})={\rm min}_{C\in\Cc^{v}}(|\widetilde{C}|)-3.$
\igno{Let $\widetilde{\Cc}=\{C \mbox{ }{\rm virtual}\mbox{ }{\rm minimal}\mbox{ }{\rm cycle}:\forall e\in E(C)\mbox{ }{\rm virtual}\mbox{ } e\in B_1(M_{F_e})\}.$
By hypothesis $\widetilde{\Cc}=\{C \mbox{ }{\rm virtual}\mbox{ }{\rm minimal}\mbox{ }{\rm cycle}\mbox{ }{\rm in}\mbox{ }G\}=\Cc$. Moreover, for all $C\in\widetilde{\Cc}$,  $\widetilde{C}$ defined as in the lemma \ref{cycle-indispensable} is a minimal cycle of $(\overline{\Gamma}^{(\Pi^*)})$. On the other hand, by theorem \ref{borne-in-max} $p_2(\B_{\overline{\Gamma}})\geq={\rm min}_{C\in\Cc}(C^{(\Pi^*)})-3$, for any $C^{(\Pi^*)}\in\Cc^{(\Pi^{*})}(C)$. Let $C_0\in\Cc$ such that for any $C_0^{(\Pi^*)}\in\Cc^{(\Pi^{*})}(C_0)$, $|C_0^{(\Pi^*)}|={\rm min}_{C\in\Cc}(C^{(\Pi^*)})$. Applying proposition \ref{borne-max} we have that $\widetilde{C_0}\in \Cc^{(\Pi^{*})}(C_0)$. Then $p_2(\B_{\overline{\Gamma}})\geq |\widetilde{C_0}|\geq {\rm min}_{C\in\widetilde{\Cc}}(|\widetilde{C}|)\leq p_2(\B_{\overline{\Gamma}})$. 
The last inequality is due to theorem \ref{superieur-borne}. By this way we have $p_2(\B_{\overline{\Gamma}})={\rm min}_{C\in\Cc}(|\widetilde{C}|)-3.$}

\qed

\begin{next}\nx{Example}\label{rest-bruns}
{\rm Let $G$ be the graph of the figure \ref{bruns-herzog2-dib} and $\Gamma=\Gamma(G)$ the clique complex generated by $G$. The subgraph $C$ of $G$ whose edge-set $E(C)=\{\{a,c\},\{c,d\},\{d,e\},\{e,a\}\}$ is the only minimal cycle of $G$. We associate to the facets $F_1=\{a,b,c\}$ and $F_2=\{c,d\}$ the following matrices:

\begin{center}
$M_{F_1}=\left(  \begin{array}{ll}
    	a&z\\
	z&c
      \end{array}\right)$ and
$M_{F_2}=\left(  \begin{array}{lll}
    	e&w&x\\
	w&x&d
      \end{array}\right)$.
\end{center}

Thus the graph $\overline{G}'$ of the figure \ref{bruns-herzog2-dib} is the $1$-skeleton of $\overline{\Gamma}^{\Pi^*}=\overline{\Gamma}'$. 
\mbox{ }\\
%\mbox{ }\\
\begin{figure}[h]  
  \centering
  \scalebox{0.70}{
  \includegraphics[viewport=230 510 300 640]{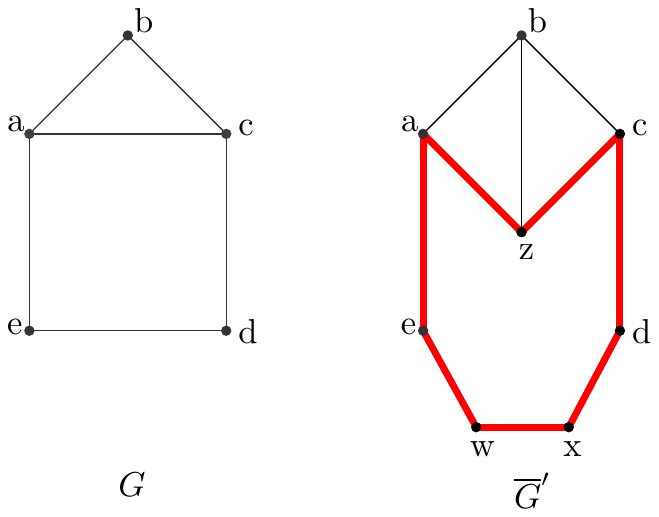}}
  \vspace{3mm}  
  \caption{} 
  \label{bruns-herzog2-dib}
\end{figure}

Besides, by a proposition of the preprint \cite{hernan}, $\J_{\overline{\Gamma}}$ is a toric ideal. Moreover, we see that the conditions of the theorem \ref{2pregul2} required to $\Gamma$ and its extension $\overline{\Gamma}$ are satisfied, so:
$$\begin{array}{ll}
	p_2(\B_{\overline{\Gamma}})&=p_2(I_{\overline{\Gamma}^{\Pi^*}})\\
	&=(|C|+|R_1^{({\Pi^*})}(C)|+\sum_{e\in R^{({\Pi^*})}_{3}(C)}\eta_e)-3\\
	&=(4+1+2)-3=4.
  \end{array}$$}
\end{next}

\end{document}